\documentclass[12pt]{article}
\usepackage{amssymb,amsmath}
\usepackage{color}

\newcommand{\R}{\mathbb{R}}
\newcommand{\N}{\mathbb{N}}
\newcommand{\F}{\mathbb{F}}

\def\half{{1\over 2}}
\def\abs#1{|#1|}

\def\norm#1{\|#1\|}
\def\cat{\mathop{\rm cat}\nolimits}
\def\dist{\mathop{\rm dist}\,}

\def\claim#1{\noindent {\bf #1}}

\def\flushright#1{{\unskip\nobreak\hfil\penalty50\hskip2em\hbox{}\nobreak\hfil%
#1\parfillskip=0pt\finalhyphendemerits=0\par}}

\def\bull{\vrule height 1.8ex width 1.0ex depth .1ex }
\def\QED{\ifmmode\eqno\hbox{$\bull$}\else\flushright{\hbox{$\bull$}}\fi}

\def\wlimit{\rightharpoonup}

\def\supp{\mathop{\rm supp}}

\def\be{\begin{equation}}
\def\ee{\end{equation}}

\def\epsilon{\varepsilon}
\def\intRN{\int_{\R^N}}

\def\Ie{I_\epsilon}
\def\Je{J_\epsilon}
\def\Qe{Q_\epsilon}

\def\calN{{\cal N}}
\def\calS{{\cal S}}

\def\calV{{\cal V}}
\def\calK{{\cal K}}

\def\whS{\widehat S}

\def\wrho{\widehat\rho}

\def\ue{u_\epsilon}
\def\pe{p_\epsilon}
\def\wtwe{\widetilde w_\epsilon}

\def\cuplength{{\rm cupl}}

\def\hdelta{\hat\delta}

\def\calXed{{\cal X}_{\epsilon,\delta}}
\def\calXp{\calXed^{E(m_0)+\hdelta}}
\def\calXm{\calXed^{E(m_0)-\hdelta}}

\newcommand{\e}{\varepsilon}

\newtheorem{theorem}{Theorem}[section]
\newtheorem{proposition}[theorem]{Proposition}

\newtheorem{lemma}[theorem]{Lemma}
\newtheorem{remark}[theorem]{Remark}
\newtheorem{definition}[theorem]{Definition}

\begin{document}

\title
{Multiplicity of positive solutions of nonlinear Schr\"odinger equations\\
concentrating at a potential well}

\author{Silvia Cingolani\\
Dipartimento di Meccanica, Matematica e Management\\
Politecnico di Bari\\
Via E. Orabona 4, 70125 Bari, Italy\\
\\
Louis Jeanjean\\
Laboratoire de Math\'ematiques (UMR CNRS 6623)\\
Universit\'e de Franche-Comt\'e\\
16, Route de Gray, 25030 Besan\c{c}on Cedex, France\\
\\
Kazunaga Tanaka\\
Department of Mathematics\\
School of Science and Engineering\\
Waseda University\\
3-4-1 Ohkubo, Shijuku-ku, Tokyo 169-8555, Japan}

\date{}

\maketitle


\begin{abstract}
We consider singularly perturbed nonlinear Schr\"odinger equations
\be \label{eq:0.1}
- \varepsilon^2 \Delta u + V(x)u = f(u), \ \  u > 0, \ \ v \in
H^1( \R^N) \ee
where $V \in C(\R^N, \R)$ and $f$ is a nonlinear term which satisfies the so-called Berestycki-Lions conditions.
We assume that there exists a bounded domain $\Omega \subset \R^N$ such that
\[ m_0 \equiv \inf_{x \in \Omega} V(x) <    \inf_{x \in \partial \Omega} V(x) \]
and we set $K = \{ x \in \Omega \ | \ V(x) = m_0\}$. For $\e >0$ small we prove the existence of at least ${\cuplength}(K) + 1$ solutions to (\ref{eq:0.1}) concentrating, as $\e \to 0$ around $K$. We remark that, under our assumptions of $f$, the search of solutions to (\ref{eq:0.1}) cannot be reduced  to the study of the critical points
of a functional restricted to a Nehari manifold.

\end{abstract}


\section{\label{section:1} Introduction}
In these last years a great deal of work has been  devoted to the study of semiclassical standing waves for
the nonlinear Schr\"odinger equation
\begin{equation}\label{seq}
i\hbar \frac{\partial \psi}{\partial t} +
\frac{\hbar^2}{2 m }\Delta \psi - V(x)\psi + f(\psi) = 0,\ \ (t,x)
\in {\R} \times \R^N
\end{equation}
where $\hbar$ denotes the Plank constant, $i$ the imaginary unit, $m$ is a positive number,
$f(\exp(i\theta)\xi) =
\exp(i\theta)f(\xi)$ for $\theta, \xi \in \R.$  A solution of the form
$\psi(x,t) = \exp(-iEt/\hbar)v(x)$ is called a standing wave.
Then, assuming $m = \frac{1}{2}$, $\psi(x,t)$ is a solution of (\ref{seq}) if and only if the
function $v$ satisfies
\be \label{eleq} {\hbar^2} \Delta v
- (V(x) - E)v + f(v) = 0 \ \ \textrm{ in }  \ \ \R^N. \ee
For small $\hbar >0$, these standing waves are referred to as semi-classical states. The limit $\hbar  \to 0$ somehow described the transition from Quantum Mechanics to Classical Mechanics.
For the physical background for this equation, we refer to the
introduction in \cite{ABC,BW1}.
\vspace{2mm}

In this paper we are interested on positive solutions of (\ref{eleq}) in $H^1(\R^N)$ for small $\hbar >0$. For simplicity and without loss
of generality, we write $V - E$ as $V$, i.e., we shift $E$ to
$0$ and we set $\hbar  = \varepsilon$. Thus, we consider the following equation
\be \label{eq:1.1}
- \varepsilon^2 \Delta v + V(x)v = f(v), \ \  v > 0, \ \ v \in
H^1( \R^N) \ee
when $\varepsilon > 0$ is  small. Throughout the paper, we assume that $N \geq 3$ and that the potential $V$ satisfies
\begin{itemize}
\item[(V1)] $ V \in C(\R^N, \R)$ and $\inf_{x \in \R^N} V(x) = \underline{V} >0$.
\end{itemize}
We observe that defining $u(x) = v(\e x)$  equation (\ref{eq:1.1}) is equivalent to \be
\label {eq:2.1} -\Delta u  +V(\e x) u = f(u), \ \ u > 0, \ \ u \in
H^1(\R^N). \ee
We shall mainly work on equation (\ref{eq:2.1}).  Note also that for each $x_0 \in \R$ and $R>0$, $V(\e x)$ converges uniformly to $V(x_0)$ on $B(x_0/ \e, R)$ as $\e \to 0$. Thus for each $x_0 \in \R^N$ we have a formal limiting problem
\be \label{eres} -\Delta U +V(x_0)U = f(U) \quad \hbox{in} \quad
\R^N, \quad   U \in H^1(\R^N). \ee
The approaches to find solutions of (\ref{eq:2.1}) when $\e > 0$ is small and to study their behavior as $\e \to 0$  can be, roughly, classified into two categories.

A first approach relies on a reduction type method. The first result in this direction was given by Floer and Weinstein \cite{FlWe} where $N=1$, $f(u)=u^3$ and $x_0$ is a non degenerate minimum or maximum of $V$. Later Oh \cite{O,O3} generalized this result to higher dimensions $N$ and $f(u)=|u|^{p-1} u$, $ 1 < p < (N+2)/(N-2)$, $N \geq 3$, $1< p < \infty$ if $N=1,2$  for non degenerate minima or maxima of $V$.
Successively in \cite{ABC},  a Liapunov-Schmidt type procedure was used to reduce $(\ref{eq:2.1})$,  for $\epsilon >0$ small, to as finite-dimensional equation that inherits the variational structure of the original problem.
Existence of semiclassical standing waves solutions were proved concentrating at  local minima or local maxima of $V$, with nondegenerate $m$ -th derivative, for some integer $m$. This result was generalized by Li \cite{Li}, where a degeneracity of any
order of the derivative is allowed.
The reduction type methods permitted to obtain other precise and striking results  in \cite{AMN,CN,CS,DLY,DKW,KW}.
However this approach relies on the uniqueness and non-degeneracy of the ground state solutions, namely of the positive least energy solutions for the autonomous problems (\ref{eres}). This uniqueness and nondegeneracy property is true for the model nonlinearity $f(\xi)= |\xi|^{p-1} \xi$, $1 < p < (N+2)/(N-2)$, $N \geq3$  and for some classes of nonlinearities, see \cite{CEF}. However, as it is shown in \cite{DDI}, it does not hold in general. For a weakening of the nondegeneracy requirement, within the frame of the reduction methods, see nevertheless \cite{D}.
\vspace{2mm}

An alternative type of approach was initiated by Rabinowitz \cite{R}. It is purely variational and do not require the nondegeneracy condition for the limit problems (\ref{eres}). In \cite{R} Rabinowitz proved, by a mountain pass argument,
the existence of positive solutions of (\ref{eq:2.1}) for small $\e
>0$ whenever
\begin{equation}\label{global}
  \liminf_{| x|\to\infty} V(x) > \inf_{x\in \R^N} V(x).
\end{equation}
These solutions concentrate around the global minimum points of
$V$ when $\e  \to 0$, as it was shown by X. Wang \cite{W}. Later,
del Pino and Felmer \cite{DF1} by introducing a penalization
approach prove a localized version of the result of Rabinowitz and
Wang (see also \cite{DF2,DF3,DF4,G} for related results). In
\cite{DF1}, assuming (V1) and,
\begin{itemize}
\item[(V2)] There is a bounded domain $\Omega \subset \R^N$ such that
\[ m_0 = \inf_{x \in \Omega} V(x) <    \inf_{x \in \partial \Omega} V(x) \]
\end{itemize}
they show the existence of a single peak solution concentrating
around the minimum points  of $V$ in $\Omega$. They assume that $f$ satisfies the assumptions
\begin{itemize}
\item[(f1)] $f\in C(\R,\R)$.
\item[(f2)] $f(0)=\lim_{\xi\to 0}{f(\xi)\over \xi}=0$.
\item[(f3)] There exists $p\in (1,{N+2\over N-2})$ such that
    $$  \lim_{\xi\to\infty}{f(\xi)\over \xi^p}=0.
    $$
\end{itemize}
In addition they require the {\it global Ambrosetti-Rabinowitz condition:}
\begin{equation} \label{super} \mbox{
for some } \mu > 2, \quad  0 < \mu\int_0^tf(s) ds < f(t)t, \quad  t  > 0
\end{equation}
and the monotonicity condition
\begin{equation}\label{mono}
        (0,\infty)\to\R; \quad \xi\mapsto {f(\xi)\over\xi}
        \quad \hbox{is strictly increasing}.
    \end{equation}
After some weakening of the conditions (\ref{super}) and (\ref{mono}) in  \cite{DF1,JT2}, it was finally shown in \cite{BJ} (see also \cite{BJT} for $N=1,2$) that (\ref{super})--(\ref{mono}) can be replaced  by

\begin{itemize}
\item[(f4)] There exists $\xi_0>0$ such that
    $$  F(\xi_0) >\half m_0\xi_0^2 \quad \mbox{where} \quad F(\xi)=\int_0^\xi f(\tau)\, d\tau.
    $$
\end{itemize}

Note that Berestycki and Lions \cite{BL} proved that there exists a least energy solution of (\ref{eres}) with $V(x_0) = m_0$ if (f1)-(f4) are satisfied . Also, using the Pohozaev's
identity, they showed that conditions (f3) and (f4) are necessary for
existence of a non-trivial solution of (\ref{eres}). Thus as far as the existence of solutions of \eqref{eq:2.1}, concentrating around local minimum points of $V$ is concerned, the results of \cite{BJ,BJT} are somehow optimal.
\vspace{2mm}

Subsequently the approach developed in \cite{BJ} was adapted in \cite{CJS} to derive a related result
for nonlinear Schr\"odinger equations with an external magnetic field. Also, very recently, it has been extended \cite{BT1, BT2, DPR} to obtain the existence of a family of solutions of \eqref{eq:2.1} concentrating, as $\e \to 0$, to a topologically non-trivial saddle point of $V$.
In that direction previous results obtained either by the reduction or by the variational ones as developed by del Pino and Felmer, can be found in \cite{DLY,DY,DF4,KW}.
\vspace{2mm}

In this work, letting $K$ being a set of local minima of $V$, we are interested in the multiplicity of positive solutions concentrating around $K$. Starting from the paper of Bahri and Coron \cite{BCo},
many papers are devoted to study the effect of the domain topology on the existence and
multiplicity of solutions
for semilinear elliptic problems.  We refer to \cite{BC1, BC2, BCP, ceramipassaseo, DW, DY1} for related studies for Dirichlet and Neumann boundary value problems.
We also refer to \cite{BW} for the study of nodal solutions.
\vspace{2mm}

For problem \eqref{eq:1.1}, the topology of the domain is trivial, but  the number of positive solutions
to \eqref{eq:1.1} is influenced by the topology of the level sets of the potential $V$.
This fact was showed by the first author and Lazzo in \cite{CL1}. Letting
$$K= \{ x \in \R^N ; \, V(x) = \inf_{x \in \R^N}V(x) \}$$
and assuming that $V \in C(\R^N, \R)$ and (\ref{global}) hold, they relate the number of positive solutions of \eqref{eq:1.1}  to the topology of the set $K$. It is assumed in \cite{CL1} that
$f(\xi)= |\xi|^{p-1} \xi$, $1 < p < (N+2)/(N-2)$ for $N>2$ and $ 1 < p <  + \infty$ if $N=2$.
In this case,  one can reduce the search of solutions of \eqref{eq:1.1}  to the existence of critical points of
a functional $\Ie$ restricted to the Nehari manifold $\calN_\epsilon$.  Here
    \begin{eqnarray*}
    \Ie(u)&=&\half\intRN \abs{\nabla u}^2+V(\epsilon x)u^2\, dx
            -\intRN F(u)\,dx, \quad F(\xi)=\int_0^\xi f(\tau)\, d\tau,\\
    \calN_\epsilon&=&\{ u\in H^1(\R^N)\setminus\{ 0\} ; \, \Ie'(u)u=0\}.
    \end{eqnarray*}
We remark that $\Ie$ is bounded from below on $\calN_\epsilon$. The multiplicity of positive solutions is obtained through a study of
the topology of the level set in $\calN_\epsilon$
    \begin{equation*}
        \Sigma_{\epsilon,h}=\{ u\in\calN_\epsilon ; \, \Ie(u)\in [c_\epsilon,c_\epsilon+h]\}
        \quad \hbox{for}\ h>0,
    \end{equation*}
where $c_\epsilon=\inf_{u\in\calN_\epsilon} \Ie(u)$.
In particular, in \cite{CL1} two maps are introduced $\phi_\varepsilon : M \to \Sigma_{\epsilon,h}$ and
$\beta : \Sigma_{\epsilon,h} \to M_\delta$ whose composition is homotopically equivalent to the embedding $j: M \to M_\delta$ for $h>0$ and $\delta> 0$ small.  Here $M_\delta$ denotes a $\delta$-neighborhood of $M$.

Hence one has
    \begin{equation*}
        \cat(\Sigma_{\epsilon,h}) \geq \cat_{M_\delta}(M)
    \end{equation*}
and this implies the existence of at least $\cat_{M_\delta}(M)$ positive solutions to $\eqref{eq:1.1}$.
Here $\cat_X (A)$ denotes the Lusternick-Schnirelmann category of $A$ in $X$ for any  topological pair $(X,A)$.
\vspace{2mm}

This approach was then extended to
nonlinear Schr\"odinger equations with competing potentials in  \cite{CL2}
and in \cite{C} for nonlinear Schr\"odinger equations with an external magnetic field.
\vspace{2mm}

In \cite{DLY} Dancer, Lam and Yan proved the existence of at least $\cat(K)$ positive solutions to $\eqref{eq:1.1}$ for $\varepsilon >0$ small, when $K$ is a connected compact local minimum or maximum set of $V$.
They also assume that
 $f(\xi)= |\xi|^{p-1} \xi$, with $1 < p < (N+2)/(N-2)$ for $N>2$ and  $ 1 < p <   \infty$  if $N=2$.
The proof of this result relying on a reduction method uses the special form of the nonlinear term $f$.
\vspace{2mm}

Successively in \cite{AMS} Ambrosetti, Malchiodi and Secchi proved if  $V$ has a
nondegenerate  (in the sense of Bott) manifold $M$  of critical points  and  $f(\xi)= |\xi|^{p-1} \xi$, with $1 < p < (N+2)/(N-2)$ for $N>2$ and  $ 1 < p <   \infty$ if $N=2$, that the equation  $\eqref{eq:1.1}$ has at least $\cuplength(M)+1$ critical points concentrating near points of $M$.   Here $\cuplength(M)$ denotes the cup-length of $M$ (see Definition $(\ref{claim:5.2})$).
If the critical points are local minima or local maxima, the above results can be sharped because $M$ does not need to be a manifold and $\cuplength(M)$ can be replaced by $\cat_{M_\delta}(M)$, for  $\delta >0$ small. The approach in \cite{AMS} relies on a perturbative variational method, which requires the uniqueness and the nondegeneracy of the limiting problem.

In this paper our aim is to study the multiplicity of positive solutions concentrating
around a set  $K$ of local minimum of $V$ under the conditions $(f1)-(f4)$.
In particular since we do not assume the monotonicity condition (\ref{mono}), we can not use a Nehari manifold.
We will introduce a method to analyze the topological difference between two level sets
of the indefinite functional $\Ie$ in a small neighborhood  of a set of expected
solutions.

\vspace{2mm}
Our main result is
\begin{theorem}\label{claim:1.1}
Suppose $N\geq 3$ and that (V1)--(V2) and (f1)--(f4) hold. Assume in addition that $\sup_{x\in\R^N}V(x)<\infty$. Then letting
$K=\{ x\in\Omega ; \, V(x)=m_0\}$, for sufficiently small $\epsilon >0$, \eqref{eq:1.1} has at least
$\cuplength(K)+1$ positive solutions $v_{\e}^i$, $i=1, \dots, \cuplength(K)+1$ concentrating  as $\e \to 0 $ in $K$,
where $\cuplength(K)$ denotes the cup-length defined with Alexander-Spanier cohomology with coefficients
in the field $\F$.

 \end{theorem}

\begin{remark}\label{examples}
If $M=S^{N-1}$, the $N-1$ dimensional sphere in $\R^N$, then $\cuplength(M) + 1=  \cat(M) =2$.  If $M=T^N$ is the $N$-dimensional torus, then $\cuplength(M) + 1 = \cat(M)= N + 1$.
However in general $\cuplength(M) +1 \leq \cat(M)$.
\end{remark}

\begin{remark}\label{mean-concen}
When we say that the solutions  $v_{\e}^i$, $i=1, \dots, \cuplength(K)+1$ of Theorem \ref{claim:1.1} concentrate when $\e \to 0$ in K, we mean that there exists a maximum point $x_{\e}^i$ of $v_{\e}^i$ such that $\lim_{\e \to 0}dist(x_{\e}^i, K) =0$ and that, for any such $x_{\e}^i$, $w_{\e}^i(x) = v_{\e}^i(\e(x+x_{\e}^i))$ converges, up to a subsequence, uniformly to a least energy solution of
$$- \Delta U + m_0 U = f(U), \quad U >0, \quad U \in H^{1}(\R^N).$$
We also have
$$v_{\e}^i(x) \leq C exp(- \frac{c}{\e}|x- x_{\e}^i|) \quad \mbox{ for some } c,\, C >0.$$
\end{remark}

\begin{remark} In addition to condition $(V1)$ the boundedness of $V$ from above is assumed in Theorem \ref{claim:1.1}. Arguing as in \cite{BJ,BT1}  we could prove Theorem \ref{claim:1.1} without this additional assumption. However, for the sake of simplicity, we assume here the boundedness of $V$.
\end{remark}

\medskip
In the proof of  Theorem \ref{claim:1.1}, two  ingredients play important roles.
The first one is a suitable choice of a neighborhood $\calXed$ of a set of expected solutions.
We remark that in a situation where the search of solutions to \eqref{eq:2.1} can be reduced to a variational problem
on a Nehari manifold, we just need to study the (global) level set
$\Sigma_{\epsilon,h}\subset\calN_\epsilon$ and we can apply a standard deformation argument
which is developed for functionals defined on Hilbert manifolds without boundary.
On the contrary, in our setting, we need to find critical points in a neighborhood
which has boundary.  Thus we need to find a neighborhood positively invariant
under a pseudo-gradient flow to develop our deformation argument.
With the aid of $\epsilon$-independent gradient estimate (Proposition \ref{claim:3.1} below),
we find such a neighborhood.
Moreover we  analyze the topological difference between two level sets
of the indefinite functional $\Ie$ in  the neighborhood $\calXed$.
To this aim  a second crucial ingredient is the definition of two maps
   \begin{eqnarray*}
    &&\Phi_\epsilon:\, ([1-s_0,1+s_0]\times K, \{ 1\pm s_0\}\times K)
        \to (\calXp, \calXm),\\
    &&\Psi_\epsilon:\, (\calXp, \calXm) \to\\
    &&\qquad ([1-s_0,1+s_0]\times\Omega([0,\nu_1]), ([1-s_0,1+s_0]\setminus\{ 1\})\times\Omega([0,\nu_1]))
    \hss
    \end{eqnarray*}
where $s_0\in (0,1)$ is small, $\Omega([0,\nu_1]) =\{ y\in\Omega;\, m_0 \leq V(y) \leq m_0 + \nu_1 \}$ for a suitable $\nu_1>0$ small, $E(m_0)$ is the least energy level associated to the limit problem
   \begin{equation}
    - \Delta u + m_0 u = f(u), \quad u \in H^1(\R^N)
   \end{equation}
and $\calXed^c=\{ u\in\calXed;\, \Je(u)\leq c\}$  for any $c\in\R$.

We show  that $\Psi_\varepsilon \circ \Phi_\epsilon$ is homotopically equivalent to the embedding $j(s,p)=(s,p)$.
We emphasize that to define such maps, center of mass and a function
$P_0$ which is inspired by the Pohozaev identity are important.

Also, differently from \cite{CL1},  it will be necessary to use the notions of category and cup-length
for maps to derive our topological result (see Remark \ref{counterexample}).

The paper is organized as follows. Section \ref{section:2} is devoted to some preliminaries. In Section \ref{section:3},  we introduce, for technical reasons, a penalized functional $J_{\e}$. In Section \ref{section:4} we define our neighborhood of expected solutions $\calXed$ and we build our deformation argument. We also introduce the two maps $\Phi_\epsilon$ and $\Psi_\epsilon$ and establish some of their important properties. Finally in Section \ref{section:5}, after having briefly recall the definitions and properties that we use of category and cup-length of a maps, we give the proof of Theorem \ref{claim:1.1} and of Remark \ref{mean-concen}.


\setcounter{equation}{0}
\section{\label{section:2}  Preliminaries}
In what follows we use the notation:
    \begin{eqnarray*}
        \norm u_{H^1} &=& \left(\intRN \abs{\nabla u}^2+u^2\, dx\right)^{1/2},\\
        \norm u_r &=& \left(\intRN \abs{u}^r\, dx\right)^{1/r} \quad
                        \hbox{for} \ r\in [1,\infty).
    \end{eqnarray*}
We study the multiplicity of solutions to \eqref{eq:2.1} via a variational method. That is, we look for
critical points of the functional $\Ie \in C^1(H^1(\R^N),\R)$ defined by
    $$  \Ie(u)=\half\norm{\nabla u}_2^2+\half\intRN V(\epsilon x)u^2\, dx
        -\intRN F(u)\, dx.
    $$
The critical points of $\Ie$ are clearly solutions of  \eqref{eq:2.1}.
Without loss of generality, we may  assume that $f(\xi)=0$ for all $\xi\leq 0$. Indeed it is then easy to see from the maximum principle that any nontrivial solution of \eqref{eq:2.1} is positive.

\medskip

\subsection{\label{subsection:2.1} Limit problems}
We introduce the notation
$$  \Omega(I) =\{ y\in\Omega;\, V(y)-m_0\in I\}
    $$
for an interval $I\subset [0,\inf_{x\in\partial\Omega}V(x)-m_0)$.
We choose a small $\nu_0>0$ such that

\begin{itemize}
\item[(i)] $0<\nu_0<\inf_{x\in\partial \Omega}V(x) -m_0$,
\item[(ii)] $F(\xi_0)>\half(m_0+\nu_0)\xi_0^2$,
\item[(iii)] $\Omega([0,\nu_0])\subset K_d$, where $d>0$ is a constant
for which Lemma \ref{claim:5.4} (Section \ref{section:5}) holds.
\end{itemize}

For any $a>0$ we also define a functional $L_a \in C^1(H^1(\R^N), \R)$ by
$$  L_a(u)=\half\norm{\nabla u}_2^2+{a\over 2}\norm u_2^2
        -\intRN F(u)\, dx.
    $$
associated to the limit problem
    \begin{equation}\label{eq:2.2}
    - \Delta u + a u = f(u), \quad u \in H^1(\R^N).
    \end{equation}
We denote by $E(a)$ the least energy level for \eqref{eq:2.2}.
That is,
    $$  E(a)=\inf\{ L_a(u);\, u\not=0,\ L_a'(u)=0\}.
    $$
In \cite{BL} it is proved that there exists a least energy solution of \eqref{eq:2.2}, for any $a >0$, if (f1)--(f4) are satisfied (here we consider (f4) with $m_0 = a$). Also it is showed that each solution of \eqref{eq:2.2} satisfies the Pohozaev's identity
\begin{equation}\label{eq:poho}
\frac{N}{N-2} \int_{\R^N} |\nabla u|^2 dx + N \int_{\R^N} a \frac{u^2}{2} - F(u) dx = 0.
\end{equation}
We note that, under our choice of $\nu_0 >0$,  $E(a)$ is attained for $a\in [m_0,m_0+\nu_0]$.  Clearly
$a\mapsto E(a);\, [m_0,m_0+\nu_0]\to \R$ is continuous and strictly increasing.
Choosing $\nu_0>0$ smaller if necessary, we may assume
    $$  E(m_0+\nu_0) < 2E(m_0).
    $$
We choose $\ell_0\in (E(m_0+\nu_0), 2E(m_0))$ and we set
    $$  S_a = \{ U\in H^1(\R^N)\setminus\{ 0\} ; \,
        L_a'(U)=0,\ L_a(U) \leq \ell_0, \ U(0)=\max_{x\in\R^N}U(x)\}.
    $$
    We also define
 $$  \whS=\bigcup_{a\in [m_0,m_0+\nu_0]} S_a.
    $$
From the proof of \cite[Proposition1]{BJ} we see that $\whS$ is compact in $H^1(\R^N)$ and that its elements have a uniform exponential
decay. Namely that there exist $C$, $c>0$ such that
\begin{equation}\label{expdecay}
      U(x) +  \, \abs{\nabla U(x)} \leq C \exp(-c\abs x) \quad \hbox{for all} \quad  U\in \whS.
\end{equation}
We also have that, for $b \in [m_0, m_0 + \nu_0],$
\begin{equation}\label{eq:2.3}
\displaystyle
\lim_{b \to m_0} \sup_{{\widetilde U} \in S_b}\inf_{U \in S_{m_0}} ||U-{\widetilde U}||_{H^1} =0.
\end{equation}
To prove (\ref{eq:2.3}) we assume by contradiction that it is false. Then there exists sequences $(b_j) \subset [m_0, m_0 + \nu_0]$
such that $b_j \to m_0$ and $(U_j)$ with $U_j \in S_{b_j}$ such that $(U_j)$ is bounded away from $S_{m_0}$.
But $(U_j) \subset  \whS$ and $\whS$ is compact. Thus, up to a subsequence, $U_j \to U_0$ strongly.
To reach a contradiction it suffices to show that $L'_{m_0}(U_0) =0$ and $L_{m_0}(U_0) \leq \ell_0$.
But since $(U_j)$ is bounded
$$0 = L'_{b_j}(U_j) = L'_{m_0}(U_j) + o(1)$$
as $b_j \to m_0$. Also, since $U_j \to U_0$, we have by continuity of $L'_{m_0}$ that
$0 = L'_{m_0}(U_j) \to L'_{m_0}(U_0)$. Thus $L'_{m_0}(U_0) =0$. Similarly we can show that
$L_{m_0}(U_0) \leq \ell_0$ and this proves that (\ref{eq:2.3}) hold.
\vspace{3mm}

In what follows, we try to find our
critical points in the following set:
    $$  \calS(r)= \{ U(x-y)+\varphi(x);\, y\in\R^N,\ U\in\whS,\ \norm\varphi_{H^1}<r\}
            \quad \hbox{for}\ r>0.
    $$

\medskip

\subsection{\label{subsection:2.2} A function $P_0$}
For a later use, we define the $C^1(H^1(\R^N), \R)$ functional
    \begin{equation}\label{eq:2.4}
        P_0(u)=\left({N(\intRN F(u)\, dx-{m_0\over 2}\norm u_2^2)
            \over{N-2\over 2}\norm{\nabla u}_2^2}\right)^{1/2}.
    \end{equation}
By the Pohozaev identity (\ref{eq:poho}), for a non trivial solution $U$ of $L_{m_0}'(U)=0$, we
have $P_0(U)=1$.  Moreover a direct calculation shows that
    $$  P_0(U(x/s))=s       \quad \hbox{for all}\ s>0.
    $$
A motivation to the introduction of $P_0$ is that it permits to estimate $L_{m_0}$ from below. The proof of the following lemma is given in \cite[Lemma 2.4]{BT1} but we recall it here for completeness.

\begin{lemma}\label{claim:2.1}
Suppose that $P_0(u)\in (0,\sqrt{N\over N-2})$.  Then we have
    $$  L_{m_0}(u) \geq g(P_0(u)) E(m_0),
    $$
where
    \begin{equation}    \label{eq:2.5}
        g(t)=\half(Nt^{N-2}-(N-2)t^N).
    \end{equation}
\end{lemma}

\medskip

\claim Proof.
First we recall, see \cite{JT1} for a proof, that $E(m_0)$ can be characterized as
    \be \label{eq:2.6}
      E(m_0)=\inf\{ L_{m_0}(u);\, u\not=0,\, P_0(u)=1\}.
    \ee
Now suppose that $s=P_0(u)\in (0,\sqrt{N\over N-2})$ and set  $v(x)=u(sx)$.  We have
$P_0(v)=s^{-1}P_0(u)=1$ and then, by \eqref{eq:2.6}, $  L_{m_0}(v)\geq E(m_0).$ This implies that $${1\over N}\intRN\abs{\nabla v}^2\, dx ={2\over N-2}\intRN (F(v)-{m_0\over 2}v^2)\, dx\geq E(m_0)$$
and thus
    \begin{eqnarray*}
    L_{m_0}(u)
    &=& {s^{N-2}\over 2}\intRN\abs{\nabla v}^2\, dx+s^N\intRN\left({m_0\over 2}v^2-F(v)\right)\,dx\\
    &=& g(s){1\over N}\intRN\abs{\nabla v}^2\, dx\\
    &\geq& g(s)E(m_0).
    \end{eqnarray*}
Since $s\in (0,\sqrt{N\over N-2})$ implies $g(s)>0$ this proves the lemma. \QED

\medskip

\begin{remark}\label{wellde}
Clearly  $P_0$ takes bounded sets to bounded sets. Thus from (\ref{eq:2.3}), choosing $\nu_0>0$ smaller if necessary, we may assume
    $$  P_0(U)\in (0,\sqrt{N\over N-2}) \quad
        \hbox{for all}\quad U\in\whS =\bigcup_{a\in [m_0,m_0+\nu_0]} S_a.
    $$
Also there exists $r_0>0$ such that
    $$  P_0(u)\in (0,\sqrt{N\over N-2}) \quad
        \hbox{for all}\quad u\in \calS(r_0).
    $$
\end{remark}
\begin{remark}\label{max}
 We remark that $g(t)$ defined in \eqref{eq:2.5} satisfies $g(t)\leq 1$ for all $ t>0$ and that $ g(t)=1$ holds if and only if $t=1.$
\end{remark}
\medskip

\subsection{\label{subsection:2.3} Center of mass in $\calS(r_0)$}

Following \cite{BT1, BT2} we introduce a center of mass in $\calS(r)$.

\begin{lemma}\label{claim:2.2}
There exists $r_0$, $R_0>0$ and a map $\Upsilon:\, \calS(r_0)\to\R^N$
such that
    $$  \abs{\Upsilon(u)-p} \leq 2R_0
    $$
for all $u(x)=U(x-p)+\varphi(x)\in \calS(r_0)$ with $p\in\R^N$, $U\in\whS$,
$\norm\varphi_{H^1}\leq r_0$.  Moreover, $\Upsilon(u)$ has the following properties
\begin{itemize}
\item[(i)] $\Upsilon(u)$ is shift equivariant, that is,
    $$  \Upsilon(u(x-y))=\Upsilon(u(x))+y   \quad \hbox{for all}\ u\in\calS(r_0)\
        \hbox{and}\ y\in\R^N.
    $$
\item[(ii)] $\Upsilon(u)$ is locally Lipschitz continuous, that is, there exists
constants $C_1$, $C_2>0$ such that
    $$  \abs{\Upsilon(u)-\Upsilon(v)} \leq C_1\norm{u-v}_{H^1}
        \quad \hbox{for all}\  u, v\in\calS(r_0) \  \hbox{with}\
        \norm{u-v}_{H^1}\leq C_2.
    $$
\end{itemize}
\end{lemma}

The proof is given in \cite{BT1, BT2} in a more complicated situation.
We give here a simple proof.

\medskip

\claim Proof.
We set $r_*=\min_{U\in\whS}\norm u_{H^1}>0$ and choose $R_0>1$ such that for $U\in\whS$
    $$  \norm U_{H^1(\abs x\leq R_0)} > {3\over 4}r_*  \quad \mbox{and} \quad
        \norm U_{H^1(\abs x\geq R_0)} < {1\over 8}r_*.
    $$
This is possible by the uniform exponential decay (\ref{expdecay}).
For $u\in H^1(\R^N)$ and $p\in\R^N$, we define
    $$  d(p,u)=\psi\left(\inf_{U\in\whS} \norm{u-U(x-p)}_{H^1(\abs{x-p}\leq R_0)}\right),
    $$
where $\psi \in C_0^\infty(\R,\R)$ is such that
    \begin{eqnarray*}
        &&\psi(r)=\begin{cases}    1   &r\in [0,{1\over 4}r_*],\\
                            0   &r\in [\half r_*,\infty),\end{cases}\\
        &&\psi(r)\in [0,1]\quad \hbox{for all}\ r\in [0,\infty).
    \end{eqnarray*}
Now let
    $$  \Upsilon(u)={\displaystyle \intRN q\, d(q,u)\, dq
            \over\displaystyle \intRN d(q,u)\, dq}
        \quad \hbox{for}\quad u\in\calS({1\over 8}r_*).
    $$
We shall show that  $\Upsilon$ has the desired property.
\vspace{2mm}

Let $u\in \calS({1\over 8}r_*)$ and write $u(x)=U(x-p)+\varphi(x)$
($p\in\R^N$, $U\in \whS$, $\norm\varphi_{H^1}\leq {1\over 8}r_*$).
Then for $\abs{q-p}\geq 2R_0$ and $\tilde U\in\whS$, we have
    \begin{eqnarray*}
        \norm{u-\tilde U(x-q)}_{H^1(\abs{x-q}\leq R_0)}
        &\geq& \norm{\tilde U(x-q)}_{H^1(\abs{x-q}\leq R_0)}\\
        &&    - \norm{U(x-p)}_{H^1(\abs{x-p}\geq R_0)}-{1\over 8}r_*\\
        &>& {3\over 4}r_*-{1\over 8}r_*-{1\over 8}r_*=\half r_*.
    \end{eqnarray*}
Thus $d(q,u)=0$ for $\abs{q-p}\geq 2R_0$.  We can also see that, for small $r>0$
    $$  d(q,u)=1 \quad \hbox{for}\ \abs{q-p}<r.
    $$
Thus $B(p,r)\subset \supp d(\cdot,u)\subset B(p,2R_0)$.
Therefore $\Upsilon(u)$ is well-defined and we have
    $$  \Upsilon(u)\in B(p,2R_0) \quad \hbox{for}\quad  u\in \calS({1\over 8}r_*).
    $$
Shift equivariance and locally Lipschitz continuity of $\Upsilon$ can be
checked easily.  Setting $r_0={1\over 8}r_*$,
we have the desired result.  \QED
\vspace{3mm}

Using this lemma we have

\begin{lemma}\label{claim:2.3}
There exists $\delta_1>0$, $r_1\in (0,r_0)$ and $\nu_1\in (0,\nu_0)$ such that
for $\epsilon>0$ small
    $$  \Ie(u) \geq E(m_0) +\delta_1
    $$
for all $u\in\calS(r_1)$ with $\epsilon\Upsilon(u)\in\Omega([\nu_1,\nu_0])$.
\end{lemma}

\claim Proof.
We set $\underline M=\inf_{U\in\whS}\norm U_2^2\in(0,\infty)$, $\overline M=\sup_{U\in\whS}\norm U_2^2
\in (0,\infty)$. The fact that $\underline M >0$ and $\overline M < \infty$ can be shown in a standard way using Pohozaev identify and assumptions (f1)--(f3).  For latter use in \eqref{eq:2.10} below, we choose $\nu_1\in (0,\nu_0)$ such that
    \be\label{eq:2.7}
        E(m_0+\nu_1)-\half(\nu_0-\nu_1)\overline M > E(m_0).
    \ee
First we claim that for some $\delta_1>0$
    \be\label{eq:2.8}  \inf_{U\in\whS} L_{m_0+\nu_1}(U) \geq E(m_0)+3\delta_1.
    \ee
Indeed, on one hand, if $U\in S_a$ with $a\in [m_0,m_0+\nu_1]$, we have
    \begin{eqnarray*}
        L_{m_0+\nu_1}(U) &=& L_a(U) +\half(m_0+\nu_1-a)\norm U_2^2\\
        &\geq& E(a) + \half(m_0+\nu_1-a)\underline M
    \end{eqnarray*}
and thus
    \be\label{eq:2.9}
        \inf_{U\in \bigcup_{a\in [m_0,m_0+\nu_1]} S_a} L_{m_0+\nu_1}(U) > E(m_0).
    \ee
On the other hand, if $U\in S_a$ with $a\in [m_0+\nu_1,m_0+\nu_0]$,
    \begin{eqnarray*}
        L_{m_0+\nu_1}(U) &=& L_a(U) +\half(m_0+\nu_1-a)\norm U_2^2\\
        &\geq& E(a) -\half(\nu_0-\nu_1)\overline M\\
        &\geq& E(m_0+\nu_1)-\half(\nu_0-\nu_1)\overline M
    \end{eqnarray*}
and using \eqref{eq:2.7}, it follows that
    \be \label{eq:2.10}
        \inf_{U \in\bigcup_{a\in [m_0+\nu_1,m_0+\nu_0]} S_a} L_{m_0+\nu_1}(U) > E(m_0).
    \ee
Choosing $\delta_1 >0$ small enough,
\eqref{eq:2.8} follows from \eqref{eq:2.9} and \eqref{eq:2.10}.
\vspace{2mm}

Now observe that, since elements in $\whS$ have uniform exponential decays,
    $$  \abs{\Ie(U(x-p))-L_{V(\epsilon p)}(U)}\to 0 \quad
        \hbox{as}\ \epsilon\to 0
    $$
uniformly in $U\in\whS$, $\epsilon p\in \Omega$.  Thus, by \eqref{eq:2.8}, for $U\in\whS$, $\epsilon p\in \Omega([\nu_1,\nu_0])$
    \begin{eqnarray}\label{lowerb}
        \Ie(U(x-p)) &=& L_{V(\epsilon p)}(U) + o(1)
                \geq L_{m_0+\nu_1}(U) + o(1) \nonumber\\
        &\geq& E(m_0) +2\delta_1 \qquad \hbox{for $\epsilon>0$ small.}
    \end{eqnarray}
If we suppose that $u(x)=U(x-p)+\varphi(x)\in\calS(r_0)$ satisfies
$\epsilon\Upsilon(u)\in \Omega([\nu_1,\nu_0])$,  then by Lemma \ref{claim:2.2}, $\epsilon p$
belongs to  a $2\epsilon R_0$-neighborhood of \\
$\Omega([\nu_1,\nu_0])$. Thus by (\ref{lowerb}) it follows that
    $$  \Ie(U(x-p)) \geq E(m_0)+{3\over 2}\delta_1 \quad \hbox{for $\epsilon >0$ small.}
    $$
Finally we observe that $I'_\epsilon$ is bounded on bounded sets uniformly in $\epsilon \in (0,1]$ and that by the
compactness of $\hat S,$ $\{ U(x-p) ;\, U\in\whS,\  p \in \R^N \}$ is  bounded in $H^1(\R^N). $ Thus choosing
 $r_1\in (0,r_0)$ small, if $u(x) = U(x-p) + \varphi (x) \in\calS(r_1)$, we have
    $$  \Ie(U(x-p)+\varphi(x)) \geq \Ie(U(x-p)) -\half\delta_1 \geq E(m_0)+\delta_1.
    $$
Thus, the conclusion of lemma holds.  \QED

\medskip

\setcounter{equation}{0}
\section{\label{section:3} A penalized functional $\Je$}
For technical reasons, we introduce a penalized functional $\Je$ following
\cite{BJ}.  We assume that $\partial\Omega$ is smooth and for $h>0$ we set
    $$  \Omega_h=\{ x\in\R^N\setminus\Omega ; \, \dist(x,\partial\Omega)<h\}\cup\Omega.
    $$
We choose a small $h_0>0$ such that
    $$  V(x)>m_0 \quad \hbox{for all}\ x\in\overline{\Omega_{2h_0}\setminus\Omega}.
    $$
Let
    $$  \Qe(u)=\left(\epsilon^{-2}\norm u_{L^2(\R^N\setminus(\Omega_{2h_0}/\epsilon))}^2
            -1\right)_+^{p+1\over 2}
    $$
and
    $$  \Je(u) = \Ie(u) +\Qe(u).
    $$

\noindent
We also define
    $$  \wrho(u)=\inf\{ \norm{u-U(x-y)}_{H^1};\, y\in\R^N,\ U\in\whS\}:\, \calS(r_0)\to\R.
    $$
As shown in Proposition \ref{claim:3.2} below, a suitable critical point of $\Je$ is a critical point of $\Ie$.
A motivation to introduce $\Je$ is a property given in Lemma \ref{claim:3.4}, which enables us to get a
useful estimate from below. \medskip

The following proposition gives a uniform estimate of $\norm{\Je'}_{H^{-1}}$ in an
annular neighborhood of a set of expected solutions, which is one of the keys of our argument.

\medskip

\begin{proposition} \label{claim:3.1}
There exists $r_2\in (0,r_1)$ with the following property:
for any $0<\rho_1<\rho_0<r_2$, there exists $\delta_2=\delta_2(\rho_0,\rho_1)>0$ such
that for $\epsilon>0$ small
    $$  \norm{\Je'(u)}_{H^{-1}} \geq \delta_2
    $$
for all $u\in\calS(r_2)$ with $\Je(u)\leq E(m_0+\nu_1)$ and $(\wrho(u),\epsilon\Upsilon(u))
\in ([0,\rho_0]\times\Omega([0,\nu_0]))\setminus([0,\rho_1]\times\Omega([0,\nu_1]))$.
\end{proposition}

\medskip

\claim Proof.
By (f1)--(f3), for any $a>0$ there exists $C_a>0$ such that
    $$  \abs{f(\xi)}\leq a\abs\xi+C_a\abs\xi^{p} \quad \hbox{for all}\ \xi\in\R.
    $$
We fix a $a_0\in (0,\half\underline V)$ and compute
    \begin{eqnarray*}
        \Ie'(u)u
        &=& \norm{\nabla u}_2^2 +\int_{\R^N} V(\epsilon x)u^2\, dx -\int_{\R^N}f(u)u\, dx \\
        &\geq& \norm{\nabla u}_2^2 +\underline V\norm u_2^2 -a_0\norm u_2^2-C_{a_0}\norm u_{p+1}^{p+1}\\
        &\geq& \norm{\nabla u}_2^2 +\half \underline V\norm u_2^2 -C_{a_0}\norm u_{p+1}^{p+1}.
    \end{eqnarray*}
Now choosing $r_2>0$ small enough there exists $c>0$ such that
    \be\label{eq:3.1}
        \norm{\nabla u}_2^2 +\half \underline V\norm u_2^2 -2^p C_{a_0}\norm u_{p+1}^{p+1}
        \geq c\norm u_{H^1}^2 \quad \hbox{for all}\quad \norm u_{H^1}\leq 2r_2.
    \ee
(For a technical reason, especially to get \eqref{eq:3.20} later, we add \lq\lq $2^p$'' in front of
$C_{a_0}$.)
In particular, we have
    \be\label{eq:3.2}
        \Ie'(u)u \geq c\norm u_{H^1}^2 \quad \hbox{for all}\quad  \norm u_{H^1}\leq 2r_2.
    \ee
Now we set
    $$  n_\epsilon =\left[{h_0\over \epsilon}\right]-1
    $$
and for each $i=1,2,\cdots,n_\epsilon$ we fix a function $\varphi_{\epsilon,i}\in C_0^\infty(\Omega)$
such that
    \begin{eqnarray*}
        &&\varphi_{\epsilon,i}(x)=\begin{cases}    1   &\hbox{if $x\in\Omega_{\epsilon,i}$,} \\
                                            0   &\hbox{if $x\not\in\Omega_{\epsilon,i+1}$,} \end{cases}\\
        &&\varphi_{\epsilon,i}(x)\in [0,1], \ \abs{\varphi'_{\epsilon,i}(x)}\leq 2
            \quad \hbox{for all}\ x\in\Omega.
    \end{eqnarray*}
Here we denote for $\epsilon>0$ and $h\in (0,2h_0/\epsilon]$
    \begin{eqnarray*}
        \Omega_{\epsilon,h} &=& (\Omega_{\epsilon h})/\epsilon \\
        &=&\{ x\in \R^N\setminus(\Omega/\epsilon);\,
        \dist(x,(\partial\Omega)/\epsilon)<h\}\cup (\Omega/\epsilon).
    \end{eqnarray*}
Now suppose that a sequence $(\ue)\subset\calS(r_2)$ satisfies
for $0<\rho_0<\rho_1<r_2$
    \begin{eqnarray}
        &&\Je(\ue) \leq E(m_0+\nu_1),                        \label{eq:3.3}\\
        &&\wrho(\ue) \in [0,\rho_0],                         \label{eq:3.4}\\
        &&\epsilon\Upsilon(\ue)\in\Omega([0,\nu_0]),         \label{eq:3.5}\\
        &&\norm{\Je'(\ue)}_{H^{-1}} \to 0.                   \label{eq:3.6}
    \end{eqnarray}
We shall prove, in several steps, that for $\epsilon>0$ small
    \be\label{eq:3.7}
        \wrho(\ue)\in [0,\rho_1]\quad \hbox{and}\quad \epsilon\Upsilon(\ue)\in \Omega([0,\nu_1]),
    \ee
from which the conclusion of Proposition \ref{claim:3.1} follows.
\vspace{2mm}

\noindent
{\sl Step 1: There exists a $i_\epsilon\in \{ 1,2,\cdots, n_\epsilon\}$ such that
    \be\label{eq:3.8}
        \norm{\ue}_{H^1(\Omega_{\epsilon,i_\epsilon+1}\setminus\Omega_{\epsilon,i_\epsilon})}^2
        \leq {4r_2^2\over n_\epsilon}.
    \ee
}

\noindent
Indeed by \eqref{eq:3.5} and the uniform exponential decay of $\whS$, we have \\
$\norm{\ue}_{H^1(\R^N\setminus(\Omega/\epsilon))}\leq 2r_2$ for $\epsilon >0$ small.  Thus
    $$  \sum_{i=1}^{n_\epsilon}\norm{\ue}_{H^1(\Omega_{\epsilon,i+1}\setminus\Omega_{\epsilon,i})}^2
        \leq \norm{\ue}_{H^1(\Omega_{\epsilon,h_0/\epsilon}\setminus(\Omega/\epsilon))}^2
        \leq 4r_2^2
    $$
and there exists $i_\epsilon\in \{ 1,2,\cdots, n_\epsilon\}$ satisfying \eqref{eq:3.8}.
\medskip

\noindent
{\sl Step 2: For the $i_\epsilon$ obtained in Step 1, we set
$$
        \ue^{(1)}(x) = \varphi_{\epsilon,i_\epsilon}(x)\ue(x) \quad \mbox{and} \quad
        \ue^{(2)}(x) = \ue(x)-\ue^{(1)}(x).  $$
Then we have, as $\epsilon\to 0$,
    \begin{eqnarray}
        &&\Ie(\ue^{(1)})=\Je(\ue)+o(1),                  \label{eq:3.9}\\
        &&\Ie'(\ue^{(1)})\to 0 \quad \hbox{in}\ H^{-1}(\R^N),  \label{eq:3.10}\\
        &&\ue^{(2)}\to 0 \quad \hbox{in}\ H^1(\R^N),     \label{eq:3.11}\\
        &&\Qe(\ue^{(2)})\to 0.                           \label{eq:3.12}
    \end{eqnarray}}

\smallskip

\noindent
Observe that
 \be\label{eq:3.13}
        \Ie(\ue) = \Ie(\ue^{(1)})+\Ie(\ue^{(2)}) +o(1).
    \ee
Indeed, by \eqref{eq:3.8}
    \begin{eqnarray*}
    &&\Ie(\ue)-(\Ie(\ue^{(1)})+\Ie(\ue^{(2)}))  \\
    &=& \int_{\Omega_{\epsilon,i_\epsilon+1}\setminus\Omega_{\epsilon,i_\epsilon}}
        \nabla(\varphi_{\epsilon,i_\epsilon}\ue)\nabla((1-\varphi_{\epsilon,i_\epsilon})\ue)
        +V(\epsilon x)\varphi_{\epsilon,i_\epsilon}(1-\varphi_{\epsilon,i_\epsilon})(\ue)^2\, dx\\
       & &-\int_{\Omega_{\epsilon,i_\epsilon+1}\setminus\Omega_{\epsilon,i_\epsilon}}
            F(\ue)-F(\ue^{(1)})-F(\ue^{(2)})\, dx \\
    &\to& 0 \qquad \hbox{as}\ \epsilon\to 0.
    \end{eqnarray*}

\noindent
Thus
    \be\label{eq:3.13'}
        \Je(\ue) = \Ie(\ue^{(1)})+\Ie(\ue^{(2)}) +\Qe(\ue^{(2)})+o(1).
    \ee
We can also see that
    \be\label{eq:3.14}
        \norm{\Ie'(\ue)-\Ie'(\ue^{(1)}) -\Ie'(\ue^{(2)})}_{H^{-1}} \to 0
        \quad \hbox{as}\ \epsilon\to 0.
    \ee
In a similar way, it follows from \eqref{eq:3.6} that, since $(\ue^{(2)})$ is bounded, that
    \be\label{eq:3.15}
        \Ie'(\ue^{(2)})\ue^{(2)} + \Qe'(\ue^{(2)})\ue^{(2)}
        = \Je'(\ue)\ue^{(2)} +o(1) = o(1).
    \ee
We note that $\norm{\ue^{(2)}}_{H^1} \leq 2r_2$ and $(p+1)\Qe(u) \leq \Qe'(u)u$ for all
$u\in H^1(\R^N)$.  Thus by \eqref{eq:3.2}
    $$  c\norm{\ue^{(2)}}_{H^1}^2 + (p+1)\Qe(\ue^{(2)})\to 0 \quad
        \hbox{as}\ \epsilon\to 0,
    $$
which implies \eqref{eq:3.11} and \eqref{eq:3.12}. Now  \eqref{eq:3.11} implies that
$\Ie(\ue^{(2)})\to 0$ and thus \eqref{eq:3.9} follows from  \eqref{eq:3.13'}.
\vspace{2mm}

Finally we show \eqref{eq:3.10}.  We choose a function $\widetilde\varphi\in C_0^\infty(\R^N)$
such that
    $$  \widetilde\varphi(x)=\begin{cases}    1   &\hbox{for $x\in \Omega_{h_0}$,}\\
                                        0   &\hbox{for $x\in \R^N\setminus\Omega_{2h_0}$.}
                            \end{cases}
    $$
Then we have, for all $w\in H^1(\R^N)$,
    \begin{eqnarray*}
    \Ie'(\ue^{(1)})w &=& \Ie'(\ue^{(1)})(\widetilde\varphi(\epsilon x)w) \\
    &=& \Ie'(\ue)(\widetilde\varphi(\epsilon x)w)
        - (\Ie'(\ue)-\Ie'(\ue^{(1)}))(\widetilde\varphi(\epsilon x)w)\\
    &=& \Je'(\ue)(\widetilde\varphi(\epsilon x)w)
        - (\Ie'(\ue)-\Ie'(\ue^{(1)}))(\widetilde\varphi(\epsilon x)w)
    \end{eqnarray*}
and it follows that
    $$  \abs{\Ie'(\ue^{(1)})w}
        \leq \norm{\Je'(\ue)}_{H^{-1}}\norm{\widetilde\varphi(\epsilon x)w}_{H^1}
        +\norm{\Ie'(\ue)-\Ie'(\ue^{(1)})}_{H^{-1}}\norm{\widetilde\varphi(\epsilon x)w}_{H^1}.
    $$
We note that by \eqref{eq:3.11} and \eqref{eq:3.14}, $\norm{\Ie'(\ue)-\Ie'(\ue^{(1)})}_{H^{-1}}
\to 0$.  Therefore, by \eqref{eq:3.6}, $\norm{\Ie'(\ue^{(1)})}_{H^{-1}}\to 0$, that is
\eqref{eq:3.10} holds true;

\medskip

\noindent
{\sl Step 3: After extracting a subsequence --- still we denoted by $\epsilon$ ---,
there exist a sequence $(\pe)\subset\R^N$ and $U\in\whS$ such that
    \begin{eqnarray}
    &&\epsilon \pe\to p_0\quad \hbox{for some}\ p_0\in \Omega([0,\nu_1]),
                                                                        \label{eq:3.16}\\
    &&\norm{\ue^{(1)}-U(x-\pe)}_{H^1}\to 0,                       \label{eq:3.17}\\
    &&\Ie(\ue^{(1)})\to L_{V(p_0)}(U) \quad \hbox{as}\ \epsilon\to 0.    \label{eq:3.18}
    \end{eqnarray}
}
\noindent
Let $q_\epsilon=\Upsilon(\ue)$.  We may assume that
    $$  \ue^{(1)}(x+q_\epsilon)\wlimit \widetilde U(x) \quad \hbox{weakly in}\ H^1(\R^N)
    $$
for some $\widetilde U\in H^1(\R^N)\setminus\{ 0\}$ and also that $\epsilon q_\epsilon\to p_0$.
From the definition of $\Upsilon$ and \eqref{eq:3.10}, it follows that
    $  L_{V(p_0)}'(\widetilde U)=0.
    $
Setting
    $$  \wtwe(x)=\ue^{(1)}(x+q_\epsilon)-\widetilde U(x)
    $$
we shall prove that $\norm{\wtwe}_{H^1}\to 0.$ We have
\begin{eqnarray}
    &&\Ie'(\ue^{(1)})\wtwe(x-q_\epsilon)\nonumber\\
    &=&\ \Ie'(\widetilde U(x-q_\epsilon)+\wtwe(x-q_\epsilon))\wtwe(x-q_\epsilon) \nonumber\\
    &=&\ \int_{\R^N} \nabla(\widetilde U+\wtwe)\nabla \wtwe
        +V(\epsilon x+\epsilon q_\epsilon)(\widetilde U+\wtwe)\wtwe\, dx \nonumber\\
    &&\    -\int_{\R^N} f(\widetilde U+\wtwe)\wtwe\,dx\nonumber\\
    &=&\ L_{V(p_0)}'(\widetilde U)\wtwe
    +\int_{\R^N} \abs{\nabla\wtwe}^2 +V(\epsilon x+\epsilon q_\epsilon)\wtwe^2\,dx\nonumber\\
    &&\ +\int_{\R^N} (V(\epsilon x+\epsilon q_\epsilon)-V(p_0))\widetilde U\wtwe\, dx
    +\int_{\R^N} (f(\widetilde U)-f(\widetilde U+\wtwe))\wtwe\,dx\nonumber\\
    &=&\  \int_{\R^N} \abs{\nabla\wtwe}^2 +V(\epsilon x+\epsilon q_\epsilon)\wtwe^2\,dx
        +(I)+(II).                              \label{eq:3.19}
    \end{eqnarray}
It is easy to see $(I)\to 0$ as $\epsilon\to 0$. Now
since $\abs{f(\xi)} \leq a_0\abs\xi +C_{a_0}\abs\xi^{p}$ for all $\xi\in\R$,
    \begin{eqnarray*}
    \abs{(II)} &\leq& \int_{\R^N} (a_0(\abs{\widetilde U}+\abs{\widetilde U+\wtwe})
        +C_{a_0}(\abs{\widetilde U}^p+\abs{\widetilde U+\wtwe}^p))
        \abs{\wtwe}\, dx\\
    &\leq& \int_{\R^N} a_0\abs{\wtwe}^2 + 2^pC_{a_0}\abs{\wtwe}^{p+1}
        +(2a_0\abs{\widetilde U}+(1+2^p)C_{a_0}\abs{\widetilde U}^p)\abs{\wtwe}
        \, dx\\
    &\leq& \int_{\R^N} a_0\abs{\wtwe}^2 + 2^pC_{a_0}\abs{\wtwe}^{p+1}
        \, dx + o(1).
    \end{eqnarray*}
Here we used the fact that $\wtwe\wlimit 0$ weakly in $H^1(\R^N)$.
Thus, by \eqref{eq:3.19} and \eqref{eq:3.10}, we have
    $$  \norm{\nabla \wtwe}_2^2 +\underline V\norm{\wtwe}_2^2
        \leq a_0 \norm{\wtwe}_2^2 + 2^p C_{a_0} \norm{\wtwe}_{p+1}^{p+1}
        +o(1)
    $$
from which we deduce, using \eqref{eq:3.1}, that
\be\label{eq:3.20}
\norm{\wtwe}_{H^1}\to 0.
\ee
At this point we have obtained  \eqref{eq:3.17}, \eqref{eq:3.18} where $\pe$ and $U$ are replaced with $q_\epsilon$ and
$\widetilde U$.
Since
\be\label{eq:3.21}
        \Ie(\ue) = \Ie(\ue^{(1)}) +o(1) = \Je(\ue) + o(1) \leq E(m_0+\nu_1)+o(1)
    \ee
implies
    $$  E(V(p_0)) \leq L_{V(p_0)}(\widetilde U) \leq E(m_0+\nu_1),
    $$
we have $p_0\in \Omega([0,\nu_1])$ and $\widetilde U$ belongs to $S_{V(p_0)}\subset\whS$ after
a suitable shift, that is, $U(x):=\widetilde U(x+y_0)\in \whS$ for some $y_0\in\R^N$.
Setting $\pe=q_\epsilon+y_0$, we get \eqref{eq:3.16}--\eqref{eq:3.18}.

\medskip

\noindent
{\sl Step 4: Conclusion}

\smallskip

\noindent
In Steps 1--3, we have shown that a sequence $(u_n) \subset \calS(r_2)$ satisfying \eqref{eq:3.3}--\eqref{eq:3.6} satisfies, up to a subsequence, and for some $U\in\whS$  \eqref{eq:3.17}--\eqref{eq:3.18} with $p_{\e} = \Upsilon(\ue) + y_0$. This implies that
    \begin{eqnarray*}
        &\epsilon\Upsilon(\ue)\to p_0\in \Omega([0,\nu_1]),\\
        &\norm{\ue(x)-U(x-\pe)}_{H^1} \to 0.
    \end{eqnarray*}
In particular since $\wrho(\ue)\to 0$ we have $\wrho(\ue) \in [0, \rho_1] $ and \eqref{eq:3.7} holds.
This ends the proof of the Proposition.  \QED

\medskip

\begin{proposition}\label{claim:3.2}
There exists $\epsilon_0>0$ such that for $\epsilon\in (0,\epsilon_0]$ if $\ue\in\calS(r_2)$
satisfies
    \begin{eqnarray}
        &&\Je'(\ue)=0,                       \label{eq:3.22}\\
        &&\Je(\ue)\leq E(m_0+\nu_1),         \label{eq:3.23}\\
        &&\epsilon\Upsilon(\ue)\in\Omega([0,\nu_0]), \label{eq:3.24}
    \end{eqnarray}
then
    \be\label{eq:3.25}
        \Qe(\ue)=0 \quad \hbox{and}\quad \Ie'(\ue)=0.
    \ee
That is, $\ue$ is a solution of \eqref{eq:1.1}.
\end{proposition}

\medskip

\claim Proof.
Suppose that $\ue$ satisfies \eqref{eq:3.22}--\eqref{eq:3.24}.  Since $\ue$ satisfies \eqref{eq:3.23} we have
    \begin{eqnarray}
        &&-\Delta \ue +
        \Bigl(V(\epsilon x)+(p+1)
        (\epsilon^{-2}\norm{\ue}_{L^2(\R^N\setminus(\Omega_{2h_0}/\epsilon))}^2-1)_+^{p-1\over 2} \nonumber\\
        &&\qquad \qquad\ \times \epsilon^{-2} \chi_{\R^N\setminus(\Omega_{2h_0}/\epsilon)}(x)\Bigr) \ue
         = f(\ue),                                  \label{eq:3.26}
    \end{eqnarray}
where $\chi_{\R^N\setminus(\Omega_{2h_0}/\epsilon))}(x)$ is the characteristic function of the set
${\R^N\setminus(\Omega_{2h_0}/\epsilon)}$.
Clearly $\ue$ satisfies \eqref{eq:3.3}--\eqref{eq:3.6} and thus, by the proof of Proposition \ref{claim:3.1}, we have
    $$  \norm{\ue}_{H^1(\R^N\setminus(\Omega_{h_0}/\epsilon))}
        \leq \norm{\ue^{(2)}}_{H^1(\R^N)}\to 0  \quad \hbox{as}\ \epsilon\to 0.
    $$
Standard elliptic estimates then shown that
    $$  \norm{\ue}_{L^\infty(\R^N\setminus(\Omega_{{3\over 2}h_0}/\epsilon))}\to 0
        \quad \hbox{as}\ \epsilon\to 0
    $$
and using a comparison principle, we deduce that for some $c$, $c'>0$
    $$  \abs{\ue(x)} \leq c'\exp(-c\dist(x,\Omega_{{3\over 2}h_0}/\epsilon)).
    $$
In particular then
    $$  \norm{\ue}_{L^2(\R^N\setminus(\Omega_{2h_0}/\epsilon))}<\epsilon \quad
        \hbox{for} \quad  \epsilon >0 \ \hbox{small}\
    $$
and we have \eqref{eq:3.25}.                     \QED
\vspace{2mm}

To find critical points of $\Je(u)$, we need the following

\medskip

\begin{proposition} \label{claim:3.3}
For any fixed $\epsilon>0$, the Palais-Smale condition holds for $\Je$ in $\{ u\in\calS(r_2);\,
\epsilon\Upsilon(u)\in\Omega([0,\nu_0])\}$.  That is, if a sequence $(u_j)\subset
H^1(\R^N)$ satisfies for some $c>0$
    \begin{eqnarray*}
        &&u_j \in\calS(r_2),\\
        &&\epsilon\Upsilon(u_j) \in \Omega([0,\nu_0]),\\
        &&\norm{\Je'(u_j)}_{H^{-1}} \to 0,\\
        &&\Je(u_j)\to c \quad \hbox{as}\ j\to \infty,
    \end{eqnarray*}
then $(u_j)$ has a strongly convergent subsequence in $H^1(\R^N)$.
\end{proposition}

\medskip

\claim Proof.
Since $\calS(r_2)$ is bounded in $H^1(\R^N)$, after extracting a subsequence if necessary, we
may assume $u_j\wlimit u_0$ weakly in $H^1(\R^N)$ for some $u_0\in H^1(\R^N)$.  We will
show that $u_j\to u_0$ strongly in $H^1(\R^N)$. Denoting $B_R=\{ x\in\R^N;\, \abs x<R\}$, it suffices to show that
    \be\label{eq:3.27}
        \lim_{R\to\infty} \lim_{j\to\infty} \norm{u_j}_{H^1(\R^N\setminus B_R)}^2=0.
    \ee
To show \eqref{eq:3.27} we first we note that, since $\varepsilon >0$ is fixed,  $\norm{u_j}_{H^1(\R^N\setminus B_L)}<2r_2$ for
a large $L>1$.  In particular, for any $n\in\N$
    $$  \sum_{i=1}^n \norm{u_j}_{H^1(A_i)}^2 < 4r_2^2,
    $$
where $A_i= B_{L+i}\setminus B_{L+i-1}$.

Thus, for any $j\in\N$, there exists $i_j\in \{1,2,\cdots,n\}$ such that
    $$  \norm{u_j}_{H^1(A_{i_j})}^2 < {4r_2^2\over n}.
    $$
Now we choose $\zeta_i\in C^1(\R,\R)$ such that $\zeta_i(r)=1$ for $r\leq L+i-1$,
$\zeta_i(r)=0$ for $r\geq L+i$ and $\zeta_i'(r)\in [-2,0]$ for all $r>0$.
We set
    $$  \widetilde u_j(x)=(1- \zeta_{i_j}(\abs x))u_j(x).
    $$
We have, for a constant  $C>0$ independent of $n$, $j$
    \begin{eqnarray}
    \Je'(u_j)\widetilde u_j &=& \Ie'(u_j)\widetilde u_j+\Qe'(u_j)\widetilde u_j, \label{eq:3.28}\\
    \Ie'(u_j)\widetilde u_j &=& \Ie'(\widetilde u_j)\widetilde u_j
        +\int_{A_{i_j}}  \nabla(\zeta_{i_j}u_j)\nabla((1-\zeta_{i_j})u_j) +V(\epsilon x)\zeta_{i_j}(1- \zeta_{i_j})u_j^2  \nonumber \\
        && +  [f((1- \zeta_{i_j})u_j) - f(u_j)](1-\zeta_{i_j})u_j\, dx \nonumber\\
    &\geq& \Ie'(\widetilde u_j)\widetilde u_j -{C\over n},       \label{eq:3.29}\\
    \Qe'(u_j)\widetilde u_j &=&
    (p+1)\left(\epsilon^{-2}\norm{u_j}_{L^2(\R^N\setminus(\Omega_{2h_0}/\epsilon))}^2-1\right)_+^{p-1\over 2}\nonumber\\
	&& \quad \times
    \int_{\R^N\setminus(\Omega_{2h_0}/\epsilon))}(1-\zeta_{i_j})u_j^2\, dx \nonumber\\
    &\geq& 0.    \label{eq:3.30}
    \end{eqnarray}
Since $\Je'(u_j)\widetilde u_j\to 0$, it follows from \eqref{eq:3.28}--\eqref{eq:3.30} that
    $$  \Ie'(\widetilde u_j)\widetilde u_j \leq {C\over n}+o(1) \quad
        \hbox{as}\ j\to \infty.
    $$
Now recording that $\norm{\widetilde u_j}_{H^1} < 2 r_2$ we have by  \eqref{eq:3.2} for some $C>0$
    $$  \norm{\widetilde u_j}_{H^1}^2 \leq {C\over n}+ o(1).
    $$
Thus, from the definition of $\widetilde u_j$, we deduce that
    $$  c\norm{u_j}_{H^1(\R^N\setminus B_{L+n})}^2
        \leq {C\over n}+o(1).
    $$
That is, \eqref{eq:3.27} holds and $(u_j)$ strongly converges.   \QED
\bigskip

The following lemma will be useful to compute the relative category.

\medskip

\begin{lemma}\label{claim:3.4}
There exists $C_0>0$ independent of $\epsilon >0$ such that
    \be\label{eq:3.31}
        \Je(u) \geq L_{m_0}(u) -C_0\epsilon^2 \quad \hbox{for all}\ u\in\calS(r_1).
    \ee
\end{lemma}

\claim Proof.
We have
    \begin{eqnarray*}
        \Je(u) &=& L_{m_0}(u)+\half\intRN(V(\epsilon x)-m_0)u^2\, dx +Q_\epsilon(u)\\
        &\geq& L_{m_0}(u) -\half(m_0-\underline V)\norm u_{L^2(\R^N\setminus(\Omega/\epsilon))}^2
                +Q_\epsilon(u).
    \end{eqnarray*}
We distinguish the two cases:  (a) $\norm u_{L^2(\R^N\setminus(\Omega/\epsilon))}^2\leq 2\epsilon^2$,
(b) $\norm u_{L^2(\R^N\setminus(\Omega/\epsilon))}^2\geq 2\epsilon^2$.
\vspace{1mm}

If case (a) occurs, we have
    $$  \Je(u) \geq L_{m_0}(u)-(m_0-\underline V)\epsilon^2
    $$
and \eqref{eq:3.31} holds.  If case (b) takes  place, we have
    $$  Q_\epsilon(u)
        \geq \left(\half\epsilon^{-2}\norm u_{L^2(\R^N\setminus(\Omega/\epsilon))}^2
                \right)^{p+1\over 2}
        \geq \half\epsilon^{-2}\norm u_{L^2(\R^N\setminus(\Omega/\epsilon))}^2
    $$
and thus
    \begin{eqnarray*}
    \Je(u) &\geq& L_{m_0}(u)
    +\half(\epsilon^{-2}-(m_0-\underline V))\norm u_{L^2(\R^N\setminus(\Omega/\epsilon))}^2\\
    &\geq& L_{m_0}(u) \quad \hbox{for $\epsilon>0$ small}.
    \end{eqnarray*}
Therefore \eqref{eq:3.31} also holds.                           \QED

\medskip

\setcounter{equation}{0}
\section{\label{section:4} A neighborhood of expected solutions}
In this section we try to find critical points of $\Je$.
First we choose a neighborhood $\calXed$ of a set of expected solutions, which is
positively invariant under a pseudo-gradient flow and in which we develop a deformation
argument.  In the sequel we will estimate a change of topology between
$\calXed\cap\{ u;\, \Je(u)\leq E(m_0)+\hdelta\}$
and $\calXed\cap\{ u;\, \Je(u)\leq E(m_0)-\hdelta\}$ using the relative category.

\medskip

\subsection{\label{subsection:4.1} A neighborhood $\calXed$}
We fix $0<\rho_1<\rho_0<r_2$ and we then choose $\delta_1$, $\delta_2>0$ according to Lemma \ref{claim:2.3}
and Proposition \ref{claim:3.1}.
We set for $\delta\in (0,\min\{ {\delta_2\over 4}(\rho_0-\rho_1),\delta_1\})$,
    $$  \calXed=\{ u\in\calS(\rho_0);\, \epsilon\Upsilon(u)\in\Omega([0,\nu_0]),\
        \Je(u)\leq E(m_0)+\delta-{\delta_2\over 2}(\wrho(u)-\rho_1)_+\}.
    $$
We shall try to find critical points of $\Je$ in $\calXed$. In this aim first note that

\begin{itemize}
\item[(a)] $u\in\calS(\rho_0)$ and $\epsilon\Upsilon(u)\in\Omega([\nu_1,\nu_0])$
imply, by Lemma \ref{claim:2.3}, that
    \be\label{eq:4.1}
        \Je(u) \geq \Ie(u) \geq E(m_0)+\delta_1\geq E(m_0)+\delta.
    \ee

\noindent
In particular,
    $$   \epsilon\Upsilon(u)\in\Omega([0,\nu_1)) \quad \hbox{for}\ u\in\calXed.
    $$
\item[(b)]  For $u\in\calXed$, if $\wrho(u)=\rho_0$, i.e., $u\in\partial\calS(\rho_0)$,
then by the choice of $\delta$
    \be\label{eq:4.2}
        \Je(u)\leq E(m_0)+\delta-{\delta_2\over 2}(\rho_0-\rho_1)\leq E(m_0)-\delta.
    \ee
\end{itemize}

\noindent
Now we consider a deformation flow defined by
    \begin{equation}\label{eq:4.3}
    \begin{cases}
        {d\eta\over d\tau}=-\phi(\eta) {\calV(\eta)\over\norm{\calV(\eta)}_{H^1}},\\
        \eta(0,u)=u,
        \end{cases}
    \end{equation}
where $\calV(u):\, \{ u\in H^1(\R^N);\, \Je'(u)\not=0\}\to H^1(\R^N)$ is a
locally Lipschitz continuous pseudo-gradient vector field satisfying
    $$  \norm{\calV(u)}_{H^1}\leq \norm{\Je'(u)}_{H^{-1}},\quad
        \Je'(u)\calV(u)\geq \half\norm{\Je'(u)}_{H^{-1}}^2
    $$
and $\phi(u):\, H^1(\R^N)\to [0,1]$ is a locally Lipschitz
continuous function. We require that $\phi(u) $ satisfies $\phi(u)=0$ if $\Je(u)\not\in
[E(m_0)-\delta,E(m_0)+\delta]$.

\medskip

\begin{proposition}\label{claim:4.1}
For any $c\in (E(m_0)-\delta,E(m_0)+\delta)$ and for any neighborhood $U$ of
$\calK_c\equiv \{ u\in \calXed;\, \Je'(u)=0,\ \Je(u)=c\}$ ($U=\emptyset$ if
$\calK_c=\emptyset$), there exist $d>0$ and a deformation
$\eta(\tau,u):\, [0,1]\times(\calXed\setminus U)\to\calXed$ such that
\begin{itemize}
\item[(i)] $\eta(0,u)=u$ for all $u$.
\item[(ii)] $\eta(\tau,u)=u$ for all $\tau\in [0,1]$ if $\Je(u)\not\in [E(m_0)-\delta,
E(m_0)+\delta]$.
\item[(iii)] $\Je(\eta(\tau,u))$ is a non-increasing function of $\tau$ for all $u$.
\item[(iv)] $\Je(\eta(1,u))\leq c-d$ for all $u\in \calXed\setminus U$ satisfying
$\Je(u)\leq c+d$.
\end{itemize}
\end{proposition}

\medskip

\claim Proof.
We consider the flow defined by \eqref{eq:4.3}. The properties (i)-(iii) follows by standard arguments from the definition \eqref{eq:4.3} and since
$\phi(u)=0$ if $\Je(u)\not\in [E(m_0)-\delta,E(m_0)+\delta]$ . Clearly also since, by Proposition \ref{claim:3.3},  $\Je$ satisfies the Palais-Smale condition for fixed $\epsilon >0$, property (iv) is standard. Thus to end the proof we just need to show that

\be\label{eq:invariant}
        \eta(\tau,\calXed)\subset\calXed \quad \hbox{for all}\ \tau\geq 0,
    \ee
namely that $\calXed $ is positively invariant under our flow. First note that because of  property  (iii), \eqref{eq:4.1} implies that for
$u\in\calXed$, $\eta(t)=\eta(t,u)$
does not enter the set $\{u;\, \epsilon\Upsilon(u)\in \Omega([\nu_1,\nu_0])\}$. Also, because of property (ii),  \eqref{eq:4.2}, show that for  $u\in\calXed$, $\eta(t)$ remains in
$\calS(\rho_0)$. Thus to show \eqref{eq:invariant} we just need to prove that the property
$$\Je(u)\leq E(m_0)+\delta-{\delta_2\over 2}(\wrho(u)-\rho_1)$$
is stable under the deformation. For this it suffices to show
that for a solution $\eta(\tau)$ of \eqref{eq:4.3},
if $0<s<t<1$ satisfies
    \begin{eqnarray*}
        &\wrho(\eta(\tau)) \in [\rho_1,\rho_0] \quad \hbox{for all}\ \tau\in [s,t],\\
        &\Je(\eta(s))\leq E(m_0)+\delta-{\delta_2\over 2}(\wrho(\eta(s))-\rho_1),
    \end{eqnarray*}
then
    $$  \Je(\eta(t))\leq E(m_0)+\delta -{\delta_2\over 2}(\wrho(\eta(t))-\rho_1).
    $$
We note that $(\wrho(\eta(\tau)), \epsilon\Upsilon(\eta(\tau)))\in [\rho_1,\rho_0]\times
\Omega([0,\nu_1])\subset ([0,\rho_0]\times\Omega([0,\nu_0]))\setminus([0,\rho_1]\times\Omega([0,\nu_1]))$
for all $\tau\in[s,t]$.
Thus by Proposition \ref{claim:3.1}, we have for $\tau\in [s,t]$
    $$  {d\over d\tau}\Je(\eta(\tau)) = \Je'(\eta){d\eta\over d\tau}
            = -\phi(\eta)\Je'(\eta){\calV(\eta)\over\norm{\calV(\eta)}_{H^1}}
            \leq -\phi(\eta){\delta_2\over 2}
    $$
    and
    \be\label{eq:4.4}
        \Je(\eta(t)) \leq \Je(\eta(s)) -{\delta_2\over 2}\int_s^t \phi(\eta(\tau))\, d\tau.
    \ee
On the other hand,
    \be\label{eq:4.5}
        \norm{\eta(t)-\eta(s)}_{H^1} \leq \int_s^t \norm{\frac{d \eta}{d \tau}}_{H^1} \, d\tau
        \leq \int_s^t \phi(\eta(\tau))\, d\tau.
    \ee
By \eqref{eq:4.4}--\eqref{eq:4.5}, and using the fact that $\abs{\wrho(\eta(t))-\wrho(\eta(s))}\leq
\norm{\eta(t)-\eta(s)}_{H^1}$, we have
    \begin{eqnarray*}
        \Je(\eta(t)) &\leq& \Je(\eta(s)) -{\delta_2\over 2}\norm{\eta(t)-\eta(s)}_{H^1}\\
        &\leq& \Je(\eta(s)) -{\delta_2\over 2}\abs{\wrho(\eta(t))-\wrho(\eta(s))}\\
        &\leq& E(m_0)+\delta-{\delta_2\over 2}(\wrho(\eta(s))-\rho_1)
                    -{\delta_2\over 2}\abs{\wrho(\eta(t))-\wrho(\eta(s))}\\
        &\leq& E(m_0)+\delta-{\delta_2\over 2}(\wrho(\eta(t))-\rho_1).
    \end{eqnarray*}
Thus $\eqref{eq:invariant}$ holds and the proof of the proposition is completed.
 \QED


\medskip

\subsection{\label{subsection:4.2} Two maps $\Phi_\epsilon$ and $\Psi_\epsilon$}
For $c\in\R$, we set
    $$  \calXed^c=\{ u\in\calXed;\, \Je(u)\leq c\}.
    $$
For $\hdelta>0$ small, using relative category, we shall estimate the change of topology between $\calXp$ and
$\calXm$.
\medskip

We recall that $K= \{ x \in \Omega ; V(x) = m_0\}. $ For $s_0\in (0,1)$ small we introduce
two maps:
    \begin{eqnarray*}
    &&\Phi_\epsilon:\, ([1-s_0,1+s_0]\times K, \{ 1\pm s_0\}\times K)
        \to (\calXp, \calXm),\\
    &&\Psi_\epsilon:\, (\calXp, \calXm) \to\\
    &&\qquad ([1-s_0,1+s_0]\times\Omega([0,\nu_1]), ([1-s_0,1+s_0]\setminus\{ 1\})\times\Omega([0,\nu_1]))
    \hss.
    \end{eqnarray*}
Here we use notation from algebraic topology: $f:\, (A,B)\to (A',B')$ means $B\subset A$,
$B'\subset A'$, $f:\, A\to A'$ is continuous and $f(B)\subset B'$.

\medskip

\noindent
{\sl Definition of $\Phi_\epsilon$}: \\
Fix a least energy solution $U_0(x)\in \whS$
of $-\Delta u+m_0 u=f(u)$ and set
    $$  \Phi_\epsilon(s,p)= U_0({x-{1\over \epsilon}p\over s}).
    $$

\noindent
{\sl Definition of $\Psi_\epsilon$}:\\
We introduce a function $\tilde P_0 : H^1(\R^N)\to \R$ by
    \begin{equation*}  \tilde P_0(u)=\begin{cases}
                            1+s_0 &\hbox{if $P_0(u)\geq 1+s_0$,}\\
                            1-s_0 &\hbox{if $P_0(u)\leq 1-s_0$,}\\
                            P_0(u)  &\hbox{otherwise,}
                        \end{cases}
    \end{equation*}
where $P_0(u)$ is given in \eqref{eq:2.4}.
\vspace{1mm}

We define our operator $\Psi_\epsilon$ by
    $$  \Psi_\epsilon(u)=(\tilde P_0(u), \epsilon\Upsilon(u))
        \quad \hbox{for}\ u\in\calXp.
    $$
In what follows, we observe that $\Phi_\epsilon$ and $\Psi_\epsilon$ are well-defined
for a suitable choices of $s_0$ and $\hdelta$.  First we deal with $\Phi_\epsilon$.
\vspace{1mm}

We fix $s_0\in (0,1)$ small so that $||U_0(\frac{x}{s}) - U_0(x)||_{H^1} < \rho_1$ for $s \in [1-s_0, 1+ s_0]$. That is for each $s \in [1- s_0, 1+ s_0]$,
$$ U_0(\frac{x}{s}) = U_0(x) + \varphi_s(x), \quad \mbox{with } ||\varphi_s||_{H^1} < \rho_1.$$
Thus
	$$U_0(\frac{x - p/\varepsilon}{s}) = U_0(x - p/\varepsilon) + \tilde\varphi_s(x) \quad \mbox{with }
	||\tilde\varphi_s||_{H^1} < \rho_1.$$
Therefore, using the first property of Lemma \ref{claim:2.2}, that is $|\Upsilon(u)-p| \leq 2R_0$ for $u(x) = U(x-p) + \varphi(x) \in \calS(\rho_0)$, we get
$$
\big|\Upsilon(U_0(\frac{x - p/\varepsilon}{s}) - p/\varepsilon \big| \leq 2 R_0.$$
It follows that, for $s\in [1-s_0,1+s_0]$
\be\label{eq:add}
        \epsilon\Upsilon(U_0({x-p/\epsilon\over s}))=p+o(1).
        \ee
Also, using Lemma \ref{claim:2.1}, we have for $p \in K$ and $s\in [1-s_0,1+s_0]$
$$ \Je(U_0({x-p/\epsilon\over s}))=L_{m_0}(U_0({x-p/\epsilon\over s}))+o(1)
        =g(s)E(m_0)+o(1).$$
Thus, choosing $\hdelta>0$ small so that $g(1\pm s_0)E(m_0)<E(m_0)-\hdelta$, which is possible by Remark \ref{max}, we see that
$\Phi_\epsilon$ is well-defined as a map
    $$  ([1-s_0,1+s_0]\times K, \{ 1\pm s_0\}\times K)
        \to (\calXp, \calXm).
    $$
Next we deal with the well-definedness of $\Psi_\epsilon$. \\
By the definition $\Psi_\epsilon$, $\Psi_\epsilon(\calXp)\subset [1-s_0,1+s_0]\times
\Omega([0,\nu_1])$.  Now we assume that $u\in\calXm$.  By Lemma \ref{claim:3.4}, we have
    $$  L_{m_0}(u) \leq \Je(u) + C_0\epsilon^2 \leq E(m_0)-\hdelta + C_0\epsilon^2.
    $$
By Lemma \ref{claim:2.1},
    $$  g(P_0(u))E(m_0) \leq L_{m_0}(u) \leq E(m_0)-\hdelta + C_0\epsilon^2.
    $$
Thus for $\epsilon >0$ small we see from Remark \ref{max} that
    $$  P_0(u) \not= 1
    $$
and $\Psi_\epsilon$ is well-defined as a map $(\calXp,\calXm)\to ([1-s_0,1+s_0]\times
\Omega([0,\nu_1]), ([1-s_0,1+s_0]\setminus\{ 1\})\times \Omega([0,\nu_1]))$.


\medskip

\noindent
The next proposition will be important to estimate $\cat(\calXp,\calXm)$.

\medskip

\begin{proposition}\label{claim:4.2}
    \begin{eqnarray*}
    &&\Psi_\epsilon\circ\Phi_\epsilon:\,([1-s_0,1+s_0]\times K, \{ 1\pm s_0\}\times K) \\
    &&\quad\to([1-s_0,1+s_0]\times\Omega([0,\nu_1]),
    ([1-s_0,1+s_0]\setminus\{ 1\})\times\Omega([0,\nu_1]))
    \end{eqnarray*}
is homotopic to the embedding $j(s,p)=(s,p)$.  That is, there exists a continuous
map
    $$  \eta:\, [0,1]\times[1-s_0,1+s_0]\times K \to [1-s_0,1+s_0]\times\Omega([0,\nu_1])
    $$
such that
    \begin{eqnarray*}
        &&\eta(0,s,p)=(\Psi_\epsilon\circ\Phi_\epsilon)(s,p), \\
        &&\eta(1,s,p)=(s,p) \quad \hbox{for all}\ (s,p)\in [1-s_0,1+s_0]\times K,\\
        &&\eta(t,s,p)\in ([1-s_0,1+s_0]\setminus\{ 1\})\times\Omega([0,\nu_1])\\
        &&\qquad \hbox{for all}\ t\in [0,1] \ \hbox{and}\ (s,p)\in \{1\pm s_0\}\times K.
    \end{eqnarray*}
\end{proposition}

\claim Proof.
By the definitions of $\Phi_\epsilon$ and $\Psi_\epsilon$, we have
    \begin{eqnarray*}
	(\Psi_\epsilon\circ\Phi_\epsilon)(s,p)
    &=&\Bigl( \tilde P_0(U_0({x-p/\epsilon\over s})),
        \epsilon\Upsilon(U_0({x-p/\epsilon\over s}))\Bigr) \\
    &=& \Bigl( s, \epsilon\Upsilon(U_0({x-p/\epsilon\over s}))\Bigr).
    \end{eqnarray*}
We set
    $$  \eta(t,s,p)=\Bigl( s, (1-t)\epsilon\Upsilon({x-p/\epsilon\over s})+tp\Bigr).
    $$
Recalling \eqref{eq:add}, we see that for  $\epsilon>0$ small $\eta(t,s,p)$ has the desired properties and
$\Psi_\epsilon\circ\Phi_\epsilon$ is homotopic to the embedding $j$.  \QED

\bigskip
\begin{remark}\label{counterexample}
As an application of the definition of category and of homotopic equivalence between maps,
one can establish that if  $X$, $\Omega^-$, $\Omega^+$ are closed sets such that
$\Omega^- \subset \Omega^+$ and $\beta : X \to \Omega^+$, $\psi : \Omega^- \to X$  are two continuous  maps
such that $\beta \circ \psi$ is homotopically equivalent to the embedding $j: \Omega^- \to \Omega^+$,
then $\cat_{\Omega^+} (\Omega^-) \leq \cat_X (X)$. See, for instance, \cite{C}.
Conversely if $X$, $X_0$, $\Omega^-$, $\Omega^-_0$, $\Omega^+$, $\Omega^+_0$ are closed sets such that
$X_0 \subset X$, $\Omega^-_0 \subset  \Omega^-$, $\Omega^+_0 \subset  \Omega^+$, $\Omega^- \subset \Omega^+$ and there exists $\Psi: (X, X_0) \to (\Omega^+, \Omega_0^+)$ and $\Phi : (\Omega^-, \Omega_0^-) \to (X, X_0)$ two continuous  maps such that $\Psi \circ \Phi$ is homotopically equivalent to the embedding $j: (\Omega^-, \Omega_0^-) \to (\Omega^+, \Omega_0^+)$, one can not infer that
$\cat_X (X, X_0) \geq  \cat_{\Omega^+}(\Omega^-,\Omega^+_0)$. Here $\cat_X(A,B)$ denotes the category of $A$ in $X$ relative to $B$, where $(X, A)$ is a topological pair and $B$ is a closed subset of $X$.

Indeed,  consider the topological pairs $(\Omega^-, \Omega_0^-) = (B, \emptyset)$,  $(\Omega^+, \Omega^+_0)= (B, S)$, $(X, X_0)=  (B, \{p\})$ where $B$ is a ball, $S = \partial B$ and $p$ is a point on the sphere $S= \partial B$.
It is possible to  construct two continuous maps $\Psi: (B,\{p\} ) \to (B, S)$ and $\Phi : (B, \emptyset) \to (B, \{p\})$ such that $\Psi \circ \Phi : (B, \emptyset) \to (B; S)$ is homotopically equivalent to the embedding
$j: (B, \emptyset) \to (B, S)$.
However $\cat_X (X,X_0)= \cat_B (B, \{p\})=0$  and $\cat_{\Omega^+} (\Omega^-,\Omega^+_0)= \cat_B(B, S) =1$.
\end{remark}

From Remark \ref{counterexample}, differently from \cite{CL1},  we can not infer, in general, that
\[
\cat(\calXp, \calXm) \geq  \cat(K, \partial K).
\]
Therefore in the work  it will be necessary to use the notions of category and cup-length for an map.

\medskip
\setcounter{equation}{0}
\section{\label{section:5} Proof of Theorem \ref{claim:1.1}}

In order to prove our theorem, we shall need some topological tools that we now present for the reader convenience.
Following \cite{BW}, see also \cite{FW1, FW2}, we define

\medskip

\begin{definition} \label{claim:5.1}
Let $B \subset A$ and $B' \subset A'$ be topological spaces and $f:(A,B) \to (A',B')$ be a continuous map, that is $f : A \to A'$ is continuous and $f(B) \subset B'$. The category $\cat(f)$ of $f$ is the least integer $k \geq 0$ such that there exist open sets $A_0$, $A_1$, $\cdots$, $A_k$ with the following properties:
\begin{itemize}
\item[(a)] $A=A_0\cup A_1\cup\cdots\cup A_k$.
\item[(b)] $B\subset A_0$ and there exists a map $h_0:\, [0,1]\times A_0\to A'$
such that
    \begin{eqnarray*}
        &&h_0(0,x)=f(x) \qquad \hbox{for all}\ x\in A_0,\\
        &&h_0(1,x)\in B' \quad \qquad \hbox{for all}\ x\in A_0,\\
        &&h_0(t,x)=f(x) \qquad \hbox{for all}\ x\in B\ \hbox{and}\ t\in [0,1].
    \end{eqnarray*}
\item[(c)] For $i=1,2,\cdots, k$, $f|_{A_i}:\, A_i\to A'$ is homotopic to a
constant map.
\end{itemize}
\end{definition}

We also introduce the cup-length of $f: (A,B) \to (A',B')$. Let $H^*$ denote Alexander-Spanier cohomology
with coefficients in the field $\F$.
We recall that the cup product $\smile$ turns $H^*(A)$ into a ring with unit $1_A$, and it turns $H^*(A,B)$
into a module over $H^*(A)$. A continuous map $f : (A,B) \to (A',B')$ induces a homomorphism
$f^* : H^*(A') \to H^*(A)$ of rings as well as a homomorphism $f^* : H^*(A',B') \to H^*(A,B)$ of
abelian groups.  We also use notation:
	$$	\tilde H^n(A')=\begin{cases} 0 &\text{for $n=0$,}\\ H^n(A')&\text{for $n>0$.}\end{cases}
	$$
For more details on algebraic topology we refer to \cite{S}.

\begin{definition} \label{claim:5.11}
For  $f:(A,B) \to (A',B')$ the cup-length, $\cuplength(f)$ is defined as follows;
when $f^*:\, H^*(A',B')\to H^*(A,B)$ is not a trivial map, $\cuplength(f)$ is defined
as the maximal integer $k\geq 0$ such
that there exist elements $\alpha_1$, $\cdots$, $\alpha_k
\in \tilde{H}^*(A')$  and $\beta\in H^*(A',B')$  with
    \begin{eqnarray*}  f^*(\alpha_1\smile\cdots\smile\alpha_k \smile \beta)
        &=& f^*(\alpha_1)\smile\cdots\smile f^*(\alpha_k)\smile f^*(\beta) \\
		&\not=& 0 \ \hbox{in}\ H^*(A,B).
    \end{eqnarray*}
When $f^*=0:\, H^*(A',B')\to H^*(A,B)$, we define $\cuplength(f)=-1$.
\end{definition}

We note that $\cuplength(f)=0$ if $f^*\not=0:\, H^*(A',B')\to H^*(A,B)$ and $\tilde{H}^*(A')=0$.

As fundamental properties of $\cat(f)$ and $\cuplength(f)$, we have

\begin{proposition}\label{claim:5.3}

\begin{itemize}
\item[(i)] For $f:\, (A,B)\to (A',B')$, $\cat(f)\geq \cuplength(f)+1$.
\item[(ii)] For $f:\, (A,B)\to (A',B')$, $f':\, (A',B')\to (A'',B'')$,
    $$  \cuplength(f'\circ f)\leq \min\{\cuplength(f'), \cuplength(f)\}.
    $$
\item[(iii)] If $f, g : (A,B) \to (A',B')$ are homotopic, then $\cuplength(f) = \cuplength(g).$
\end{itemize}
\end{proposition}

\medskip

The proof of these statements can be found in  \cite[Lemma 2.7]{BW}, \cite[Lemma 2.6 (a)]{BW} and \cite[Lemma 2.6(b)]{BW}  respectively.
Finally we recall

\begin{definition} \label{claim:5.2}
For a set $(A,B)$, we define the relative category $\cat(A,B)$ and the relative
cup-length $\cuplength(A,B)$ by
    \begin{eqnarray*}
        &&\cat(A,B)=\cat(id_{(A,B)}:\, (A,B)\to (A,B)),\\
        &&\cuplength(A,B)=\cuplength(id_{(A,B)}:\, (A,B)\to (A,B)).
    \end{eqnarray*}
We also set
    $$  \cat(A)=\cat(A,\emptyset),\quad \cuplength(A)=\cuplength(A,\emptyset).
    $$
\end{definition}

\medskip

The following lemma which is due to Bartsch \cite{B} is one of the keys of our proof and
we make use of the continuity property of Alexander-Spanier cohomology.

\medskip

\begin{lemma} \label{claim:5.4}
Let $K\subset \R^N$ be a compact set.  For a $d$-neighborhood $K_d=\{ x\in \R^N;\,
\dist(x,K)\leq d\}$ and $I=[0,1]$, $\partial I=\{0,1\}$, we consider the inclusion
    $$  j:\, (I\times K,\partial I\times K)\to (I\times K_d,\partial I\times K_d)
    $$
defined by $j(s,x)=(s,x)$.  Then for $d>0$ small,
    $$  \cuplength(j) \geq \cuplength(K).
    $$
\end{lemma}

\claim Proof.
Let $k=\cuplength(K)$ and let $\alpha_1$, $\cdots$, $\alpha_k\in H^*(K)$ ($*\geq 1$) be such
that $\alpha_1\cup\cdots\cup\alpha_k\not= 0$.  By the continuity property of Alexander-Spanier
cohomology (see \cite[Theorem 6.6.2]{S}), for $d>0$ small there exists $\alpha_1^d$, $\cdots$, $\alpha_k^d\in H^*(K_d)$
such that $i_d^*(\alpha_i^d)=\alpha_i$ for $i=1,2,\cdots,k$, where $i_d:\, K\to K_d$ is the
inclusion.

By the K\"unneth formula, the cross products give us the following isomorphisms:
    \begin{eqnarray*}
        &&\times:\, H^0(I)\otimes H^n(K) \simeq H^n(I\times K), \\
        &&\times:\, H^1(I,\partial I)\otimes H^n(K)\simeq H^{n+1}(I\times K, \partial I\times K),\\
        &&\times:\, H^0(I)\otimes H^n(K_d) \simeq H^n(I\times K_d), \\
        &&\times:\, H^1(I,\partial I)\otimes H^n(K_d)\simeq H^{n+1}(I\times K_d, \partial I\times K_d).
    \end{eqnarray*}
Let $\tau\in H^1(I,\partial I)\simeq \F$ be a non-trivial element and for any set $A$
we denote by $1_A\in H^0(A)$ the unit of the cohomology ring $H^*(A)$.
Set $\tilde \beta=\tau\times 1_{K_d}\in H^1(I\times K_d,\partial I\times K_d)$ and
$\tilde \alpha_i=1_I\times\alpha_i^d\in H^*(I\times K_d)$, then we have
    $$  j^*(\tilde \beta)=\tau\times 1_K, \quad j^*(\tilde \alpha_i)=1_I\times \alpha_i.
    $$
Thus
    \begin{eqnarray*}
        j^*(\tilde\beta\smile\tilde\alpha_1^d\smile\cdots\smile\alpha_k^d)
        &=& j^*(\tilde\beta)\smile j^*(\tilde\alpha_1^d)\smile \cdots\smile j^*(\alpha_k^d) \\
        &=& (\tau\times 1_K)\cup(1_I\times\alpha_1)\smile\cdots\smile(1_I\times \alpha_k) \\
        &=&\tau\times(\alpha_1\smile\cdots\smile\alpha_k)\not= 0.
    \end{eqnarray*}
Thus we have $\cuplength(j) \geq \cuplength(K)$.  \QED
\vspace{3mm}

Now we have all the ingredients to give the
\vspace{2mm}

\claim Proof of Theorem \ref{claim:1.1}.
We observe that for $\epsilon>0$ small
    \be\label{eq:5.1}
        \#\{u\in \calXp\setminus\calXm;\, \Je'(u)=0\}
        \geq \cat(\calXp,\calXm).
    \ee
Using Proposition \ref{claim:4.1}, \eqref{eq:5.1} can be proved in a standard way (c.f. Theorem {4.2} of
\cite{FW1}).
\vspace{2mm}

By (i) of Proposition \ref{claim:5.3}, we have
    \be\label{eq:5.2}
        \cat(\calXp,\calXm) \geq \cuplength(\calXp,\calXm)+1.
    \ee
Since $\Psi_\epsilon\circ\Phi_\epsilon=\Psi_\epsilon\circ id_{(\calXp,\calXm)}\circ\Phi_\epsilon$,
it follows from (ii) of Proposition \ref{claim:5.3} that
    \begin{eqnarray}
        \cuplength(\Psi_\epsilon\circ\Phi_\epsilon)
        &\leq& \cuplength(id_{(\calXp,\calXm)}) \nonumber\\
        &=& \cuplength(\calXp,\calXm). \label{eq:5.3}
    \end{eqnarray}
Now recall that, by our choice of $\nu_0$, we have $\Omega([0,\nu_1])\subset K_d$.  Thus, letting
\begin{eqnarray*}
    &\tau:\,&([1-s_0,1+s_0]\times \Omega([0,\nu_1]), ([1-s_0,1+s_0]\setminus\{ 1\})\times\Omega([0,\nu_1])) \\
    &&\to([1-s_0,1+s_0]\times K_d,
    ([1-s_0,1+s_0]\setminus\{ 1\})\times K_d)
    \end{eqnarray*}
    be the inclusion we see, using Proposition  \ref{claim:4.2}, that $\tau \circ \Psi_\epsilon\circ\Phi_\epsilon$ is homotopic to the inclusion
 $ j:\, (I\times K,\partial I\times K)\to (I\times K_d,\partial I\times K_d)$. Thus on one hand, by (ii) of Proposition  \ref{claim:5.3},
 \be\label{eq:5.44}
        \cuplength(\Psi_\epsilon\circ\Phi_\epsilon) \geq \cuplength( \tau \circ \Psi_\epsilon\circ\Phi_\epsilon).
    \ee
    On the other hand, by (iii) of Proposition  \ref{claim:5.3},
    \be\label{eq:5.4}
        \cuplength(\tau \circ \Psi_\epsilon\circ\Phi_\epsilon) = \cuplength(j).
    \ee
At this point using \eqref{eq:5.2}--\eqref{eq:5.4} and Lemma \ref{claim:5.4},
we deduce that
    $$  \cat(\calXp,\calXm) \geq \cuplength(K)+1.
    $$
Thus by \eqref{eq:5.1}, $\Je$ has at least $\cuplength(K)+1$ critical points in $\calXp\setminus\calXm$.
Recalling Proposition \ref{claim:3.2}, this completes the proof of the Theorem.  \QED

\medskip

\claim Proof of Remark \ref{mean-concen}.
From the proof of Proposition \ref{claim:3.1} we know that for any $\nu_0 >0$ small enough  the critical points $u_{\e}^i$, $i= 1, \dots , \cuplength(K)+1$ satisfy
$$||u_{\e}^i (x) - U^i (x- x_\e^i) ||_{H^1} \to 0$$
where $\e x_\e^i = \e \Upsilon(u_{\e}^i) + o(1) \to x_0^i \in ([0, \nu_0])$ and $U^i \in \whS$. Thus $w_{\e}^i(x) = u_{\e}^i(x + x_\e^i)$ converges to $U^i \in \whS$. Now observing that these results holds for any $\nu_0 >0$ and any $\ell_0 > E(m_0+ \nu_0)$ we deduce, considering sequences $\nu_0^n \to 0$, $\ell_0^n \to E(m_0)$ and making a diagonal process, that it is possible to assume that each $w_{\e}^i$ converges to a least energy solution of
$$- \Delta U + m_0 U = f(U), \quad U >0, \quad U \in H^1(\R^N).$$
Clearly also $$u_{\e}^i(x) \leq C exp(- c|x - x_{\e}^i|), \quad \mbox{for some } c, C >0$$
and this ends the proof. \QED
\vspace{2mm}

\noindent
{\bf Acknowledgments.}
The authors would like to thank Professor Thomas Bartsch for providing them with a proof of Lemma \ref{claim:5.4} and for very helpful discussions on relative category and cup-length.
The first author would like to thank Professor Marco Degiovanni for very helpful discussions on relative category and for suggesting the counter-example in Remark \ref{counterexample}.


\bigskip

\begin{thebibliography} {l}
\footnotesize

\bibitem{ABC} A. Ambrosetti, M. Badiale, S. Cingolani,
Semiclassical states of nonlinear Schr\"odinger equations,
Arch. Rational Mech. Anal. 140 (1997), no. 3, 285--300.

\bibitem{AMN} A. Ambrosetti, A. Malchiodi, W.M. Ni., Singularly perturbed elliptic equations with symmetry: existence of solutions concentrating on spheres I, Comm. Math. Phys.,
235 (2003), 427-466.


\bibitem{AMS} A. Ambrosetti, A. Malchiodi, S. Secchi, Multiplicity
results for some nonlinear {S}chr\"{o}dinger equations with
potentials,  Arch. Ration. Mech. Anal.  159 (2001), 253--271.


\bibitem{B} T. Bartsch, Note on category and cup-length, private communication 2012.

\bibitem{BL} H. Berestycki, P.-L. Lions,
Nonlinear scalar field equations. I. Existence of a ground state,
Arch. Rational Mech. Anal. 82 (1983), no. 4, 313--345.

\bibitem{BCo} A. Bahri, J. M. Coron,
On a nonlinear elliptic equation involving the critical Sobolev exponent: the effect of the topology of the domain,
Comm. Pure Appl. Math. 41 (1988), no. 3, 253--294.


\bibitem{BW} T. Bartsch, T. Weth,
The effect of the domain's configuration space on the number of nodal solutions
of singularly perturbed elliptic equations,
Topol. Methods Nonlinear Anal. 26 (2005), 109--133.

\bibitem{BC1} V. Benci, G. Cerami,
The effect of the domain topology on the number of positive solutions of
nonlinear elliptic problems,
Arch. Rational Mech. Anal. 114 (1991), no. 1, 79--93.

\bibitem{BC2} V. Benci, G. Cerami,
Multiple positive solutions of some elliptic problems via the Morse theory and
the domain topology,
Calc. Var. Partial Differential Equations 2 (1994), no. 1, 29--48.

\bibitem{BCP} V. Benci, G. Cerami and D. Passaseo, On the number of the positive solutions of some nonlinear elliptic problems, Nonlinear Analysis, tribute in honor of G. Prodi, Quaderno Scuola Normale Sup. Pisa 1991, pp.93-107.


\bibitem{BJ} J. Byeon, L. Jeanjean,
Standing waves for nonlinear Schr\"odinger equations with a general nonlinearity,
Arch. Ration. Mech. Anal. 185 (2007), no. 2, 185--200.

\bibitem{BJT} J. Byeon, L. Jeanjean, K. Tanaka,
Standing waves for nonlinear Schr\"odinger equations with a general nonlinearity:
one and two dimensional cases,
Comm. Partial Differential Equations 33 (2008), Issue 6, 1113--1136.

\bibitem{BT1} J. Byeon, K. Tanaka,
Semi-classical standing waves for nonlinear Schr\"odinger equations
at structurally stable critical points of the potential,
J. Euro Math. Soc. to appear.

\bibitem{BT2} J. Byeon, K. Tanaka,
Semiclassical standing waves with clustering peaks for nonlinear Schr\"odinger equations,
Memoir Amer. Math. Soc., to appear.

\bibitem{BW1} J. Byeon, Z.-Q. Wang, Standing waves with a critical
frequency for nonlinear Schr\"odinger equations, Arch. Ration.
Mech. Anal.,  165, (2002), 295--316.

%
%



\bibitem{ceramipassaseo} G. Cerami, D. Passaseo, Existence and multiplicity of positive solutions
for nonlinear elliptic problems in exterior domains with "rich" topology, Nonlin. Anal. TMA 18 (1992), 109--119.



\bibitem{C} S. Cingolani,
Semiclassical stationary states of nonlinear Schr\"odinger equations with
an external magnetic field, J. Differential Equations 188 (2003), no. 1, 52--79.



\bibitem{CJS} S. Cingolani, L. Jeanjean, S. Secchi,
Multi-peak solutions for magnetic {NLS} equations without non-degeneracy conditions,
ESAIM Control Optim. Calc. Var. 15 (2009), no. 3, 653--675.

\bibitem{CL1} S. Cingolani, M. Lazzo,
Multiple semiclassical standing waves for a class of nonlinear Schr\"odinger equations,
Topol. Methods Nonlinear Anal. 10 (1997), no. 1, 1--13.

\bibitem{CL2} S. Cingolani, M. Lazzo,
Multiple positive solutions to nonlinear Schr\"odinger equations with competing potential functions,
J. Differential Equations  160 (2000), no. 1, 118--138.


\bibitem{CN} S. Cingolani, M. Nolasco, Multi-peak periodic semiclassical states for a class of nonlinear Schr\"dinger equations, Proceedings of the Royal Society of Edinburgh, vol. 128 A (1998),  1249-1260.


\bibitem{CS} S. Cingolani, S. Secchi, Semiclassical limit for nonlinear Schr\"odinger equations with magnetic fields, Journal of Mathematical Analysis and Applications, vol. 275 (2002),  108-130



\bibitem{CEF} C. Cortazar, M. Elgueta, P. Felmer, Uniqueness of positive solutions of $\Delta u + f(u) = 0$ in $\R^N$, $N \geq 3$,
Arch. Ration. Mech. Anal.,  142, (1998), 127--141.

\bibitem{D} E. N. Dancer, Peak solutions without non-degeneracy conditions, J. Diff. Equa., {\bf 246}, (2009), 3077-3088.

\bibitem{DLY} E. N. Dancer, K. Y. Lam, S. Yan
The effect of the graph topology on the existence of multipeak solutions for nonlinear
Schr\"odinger equations, Abstr. Appl. Anal. 3 (1998), no. 3-4, 293--318.

\bibitem{DW} E. N. Dancer, J. Wei,
On the effect of domain topology in a singular perturbation problem,
Top. Methods Nonlinear Analysis, 4 (1999), 347--368.

\bibitem{DY1} E. N. Dancer, S. Yan,
A singularly perturbed elliptic problem in bounded domains with nontrivial topology,
Adv. Differential Equations 4 (1999), no. 3, 347--368.



\bibitem{DY} E. N. Dancer, S. Yan, On the existence of multipeak
solutions for nonlinear field equations on $\R^N$, Discrete
Contin. Dynam. Systems,  6, (2000), 39--50.

\bibitem{DDI} J. Davila, M. del Pino, I. Guerra, Non-uniqueness of positive ground states of non-linear
Schr\"odinger equations, Proc. London Math. Soc. (3) 106 (2013) 318-344.


\bibitem{DPR} P. D'Avenia, A. Pomponio, D. Ruiz, Semiclassical states for the nonlinear Schr\"odinger equation on saddle points of the potential via variational methods, J. Funct. Anal. 262 (2012), 4600--4633.

\bibitem{DF1} M. del Pino, P. Felmer,
Local mountain passes for semilinear elliptic problems in unbounded domains,
Calc. Var. Partial Differential Equations 4 (1996), no. 2, 121--137.

\bibitem{DF2} M. del Pino, P. Felmer,
Semi-classical states for nonlinear {S}chr\"{o}dinger equations, J.
Funct. Anal. 149 (1997), 245--265.

\bibitem{DF3} M. del Pino, P. L. Felmer, Multi-peak bound states for
nonlinear Schr\"odinger equations, Ann. Inst. Henri Poincar\'e,
 15, (1998), 127--149.

\bibitem{DF4} M. del Pino, P. L. Felmer, Semi-classical states  for
nonlinear Schr\"odinger  equations: a variational reduction
method, Math. Ann., 324, 1, (2002), 1-32.

\bibitem{DKW} M. del Pino, M. Kowalczyk, J. Wei, Concentration on curves for nonlinear Schr\"odinger equations, Comm. Pure Appl. Math, {\bf 60}, (2007), 113-146.

\bibitem{FlWe} A. Floer, A. Weinstein,
Nonspreading wave packets for the cubic Schr\"odinger equation with a bounded potential,
J. Funct. Anal., 69 ,(1986), no. 3, 397--408.


\bibitem{FW1}
G. Fournier, M. Willem,
Multiple solutions of the forced double pendulum equation,
Ann. Inst. H. Poincar\'e Anal. Non Lineaire 6 (1989), suppl., 259--281.

\bibitem{FW2}
G. Fournier, M. Willem,
Relative category and the calculus of variations, Variational methods,
Progr. Nonlinear Differential Equations Appl. Birkhauser, Boston 4, 95--104 (1990).

\bibitem{G} C. Gui, Existence of multi-bump solutions for nonlinear
Schr\"odinger equations via variational method, Comm. in
P.D.E., 21, (1996), 787--820.


\bibitem{KW} X. Kang and J. Wei, On interacting bumps of semi-classical
states of nonlinear Schr\"odinger equations, Adv. Differential
Equations, 5, (2000), 899--928.

\bibitem{JT1} L. Jeanjean, K. Tanaka,
A remark on least energy solutions in $\R^N$,
Proc. Amer. Math. Soc. 131, 8, (2003) 2399--2408.

\bibitem{JT2} L. Jeanjean, K. Tanaka,
Singularly perturbed elliptic problems with superlinear or asymptotically linear nonlinearities,
Calc. Var. Partial Diff. Eq. 21, No. 3, 287--318 (2004).


\bibitem{Li} Y. Y. Li, On a singularly perturbed elliptic
    equation, Adv.  Diff. Equations 2 (1997), 955--980.

\bibitem{O} Y. G. Oh, Existence of semiclassical bound states of
    nonlinear {S}chr\"{o}dinger equations, Comm. Partial Diff. Eq.
  13 (1988), 1499--1519.

\bibitem{O3} Y.G. Oh, On positive multi-lump bound states of nonlinear Schr\"odinger
equations under multiple well potential, Comm. Math. Phys., 131, (1990), 223--253.

\bibitem{R} P. H. Rabinowitz,
On a class of nonlinear Schr\"odinger equations, Z. Angew. Math. Phys. 43 (1992), no. 2, 270--291.

\bibitem{S} E. H. Spanier,  Algebraic topology, McGraw Hill, New York, (1966).




\bibitem{W} X. Wang, On concentration of positive bound states of
nonlinear Schr\"odinger equations, Comm. Math. Phys., 153,
(1993), 229--244.


\end{thebibliography}
\end{document}